\pdfoutput=1
\documentclass[10pt]{amsart}
\usepackage{latexsym}
\usepackage{amsmath, accents}
\usepackage{amssymb, amscd}
\usepackage{ifthen}
\usepackage{amsthm}
\input epsf

\newtheorem{dfn}{Definition}[section]

\newtheorem{rems}[dfn]{Remarks}
\newtheorem{thm}{Theorem}

\newtheorem{theor}[dfn]{Theorem}

\newtheorem{lem}[dfn]{Lemma}

\newtheorem{prop}[dfn]{Proposition}
\newtheorem{cor}[dfn]{Corollary}

\newcommand\opn{\mathrel{\ooalign{$\subseteq$\cr\hidewidth\raise.225ex\hbox{$\circ\mkern.5mu$}\cr}}}

\def\proof{\par\medskip\noindent{\it Proof: }}

\def\Rc{{\mathcal R}}
\def\H{{\Bbb H}}
\def\Z{{\Bbb Z}}

\def\Q{{\mathcal Q}}

\def\W{{\mathcal W}}

\def\N{{\Bbb N}}

\def\F{{{\mathcal F}}}
\def\T{{{\mathcal T}}}
\def\P{{{\mathcal P}}}
\def\U{{{\mathcal U}}}
\def\W{{{\mathcal W}}}
\def\V{{{\mathcal V}}}
\def\c{\overline{\C}}

\def\al{\alpha}

\def\ga{\gamma}
\def\g{\gamma}
\def\Th{\Theta}
\def\g{\gamma}

\def\G{{\Gamma}}

\def\d{\delta}

\def\De{\Delta}

\def\si{\sigma}

\def\C{{\mathcal C}}

\def\la{\lambda}
\def\La{{\bf\Lambda}}

\def\si{\sigma}

\def\d{\delta}

\def\ti{\widetilde}

\def\b{\mathsf{(b)}}

\def\iso4{{\rm Isom}{\H^4}}

\def\h3{{{\Bbb H}^3}}
\def\h4{{\H^4}}

\def\sp{\S^{n-1}_\infty}
\def\sp3{{\S}^{3}_\infty}

\def\bx{$\hfill\square$}
\def\ent{\rm Ent}

\def\d{\delta}

\def\bw{{\bf w}}

\def\bu{{\mathbf u}}
\def\bv{{\bf v}}

\def\act{\curvearrowright}

\def\b{\mathbf}
\def\s{\mathsf}
\def\d{\boldsymbol\delta}
\def\ti{\widetilde}

\def\iso3{{\rm Isom}{\H^3}}

\def\bsn{{\mathsf S}^nX}
\def\bs2{{\mathsf S}^2X}

\def\bde2{{\bf \Delta}^2X}

\def\T{\mathcal T}
\def\eG{{\rm EG}(\Gamma)}
\def\HL{{\rm Hull}(\La H)}
\def\Par{{\rm Par}}
\let\s\mathsf
\def\foral{\forall\hspace*{0.5mm}}
\def\exist{\exists\hspace*{0.7mm}}

\let\s\mathsf
\let\b\mathbf
\let\c\mathcal
\def\X{\widetilde X}

\subjclass[2010]{Primary 20F65, 20F67; Secondary 57M07, 22D05}

\title{Similar relatively hyperbolic actions of a group}

\author{Victor Gerasimov}
\address{Victor Guerassimov, Departamento de Matem\'atica,
Universidade Federal de Minas Gerais, Av. Antonio Carlos, 6627/
Caixa Postal 702, CEP 31270-901 Belo Horizonte, MG, Brasil}
\email{victor@mat.ufmg.br}
\author{Leonid Potyagailo}
\address{ Leonid Potyagailo, UFR de Math\'ematiques, Universit\'e
de Lille 1, 59655 Villeneuve d'Ascq cedex, France}
\email{potyag@math.univ-lille1.fr}

\begin{document}

\def\pdFile#1#2#3{\ifnum\pdfoutput>0%
\pdfximage width #2pt height #3pt{#1.pdf}\pdfrefximage\pdflastximage\else\kern#2pt\vbox to #3pt{\vss}\fi}

\date{\today}

\begin{abstract}
Let  a discrete group $G$ possess two convergence actions by
homeomorphisms on compacta $X$ and $Y$. Consider the following
question: does there exist a convergence action
$G{\curvearrowright}Z$ on a compactum $Z$ and continuous
equivariant maps $X\leftarrow Z\to Y$? We call the space $Z$ (and
action of $G$ on it) {\it pullback} space (action). In such
general setting a negative answer follows from a recent result of
O.~Baker and T.~Riley [BR].

Suppose, in addition, that the initial actions are relatively
hyperbolic that is they are non-parabolic and the induced action
on  the distinct pairs are cocompact. Then the existence of the
pullback space if $G$ is finitely generated follows from
\cite{Ge2}. The main result of the paper claims that the pullback
space exists if and only if  the maximal parabolic subgroups of
one of the actions  are dynamically quasiconvex for the other one.

 We provide an example
of two relatively hyperbolic actions of the free group $G$ of
countable rank for which the pullback action does not exist.

We study    an analog of the notion of geodesic flow for
relatively hyperbolic groups. Further these results are used to
prove the main theorem.

\end{abstract}

\maketitle

\bigskip

\section{Introduction}

This paper is a further development of our project of studying
convergence group actions including the actions of relatively
hyperbolic groups.

 An action  of a discrete group $G$ by homeomorphisms of a
compactum $X$ is
 said to
\it have convergence property \rm if the induced action on the
space of distinct triples of $X$ is properly discontinuous. We
call such an action 3\it-discontinuous\rm. The complement
$\bold\Lambda_XG$ of the maximal open subset where the action is
properly discontinuous is called the \it limit set \rm of the
action. The action is said to be \it minimal \rm if
$\bold\Lambda_XG{=}X$.

  The goal of the paper is to establish similarity properties
between different convergence actions of a fixed group. The first
motivation for us was  the following question: \vskip3pt $\s {\sf
Q1}:$ \sl Given two minimal $3$-discontinuous actions of a group
$G$ on   compacta  $X,Y$ does  there exists a $3$-discontinuous
action $G\act Z$ on a compactum $Z$ and continuous equivariant
maps $X\leftarrow Z\to Y$?\rm \vskip3pt We call such an action \it
pullback \rm action and the space $Z$ \it pullback \rm space.

 The answer to this question is negative in
general. This follows from a recent result of  O.~Baker and
T.~Riley \cite{BR}. They indicated a hyperbolic group $G$ and a
free subgroup $H$
  of rank three such that the embedding $H\to G$ does not admit
  an equivariant continuous extension to the  hyperbolic boundaries
  $\partial_\infty H\to\partial_\infty G$ (so called {\it Cannon-Thurston
  map}).
  It is an easy consequence of this result  that $3$-discontinuous actions $H{\curvearrowright}\partial_\infty H$
and $H{\curvearrowright}\bold\Lambda_{\partial_\infty G}H$ do not
possess a pullback 3-discontinuous action  (see section 4).

\vskip3pt   The question $\s {\sf Q1}$ has a natural modification.
Suppose in addition  that our actions $G{\curvearrowright}X$,
$G{\curvearrowright}Y$ are $2$\it-cocompact\rm,  that is the
quotient of the space of {\it distinct pairs} of the corresponding
space by $G$  is compact. This condition is natural because the
class of groups which admit non-trivial 3-discontinuous and
2-cocompact actions (we say $32$\it-actions\rm) coincides with the
class of {\it relatively hyperbolic} groups \cite[Theorem
3.1]{GePo2}.

So the following question is the main subject of the paper.
\vskip3pt $\s Q2:$ \sl Given two minimal $32$-actions of a group
$G$ on compacta $X,Y$ does there exist a $3$-discontinuous action
$G{\curvearrowright}Z$ possessing continuous equivariant maps
$X\leftarrow Z\to Y$?\rm \vskip3pt

  We note that if such an action $G\act Z$ exists then one can choose it
  to be 2-cocompact (see  Lemma \ref{quotient}). So  the
pullback action is also  of type $(32)$.

Recall few standard
  definitions.
 An action on a compactum is called {\it parabolic} if it admits
  a unique fixed point. For a 3-discontinuous action $G\act X$ a
   point $p\in X$ is called {\it
  parabolic} if it is the unique fixed point for its stabilizer
  $P={\rm St}_G p$ and ${\rm St}_G p$ is infinite.
  The subgroup $P$ is  called a  {\it maximal parabolic} subgroup.
 We denote by $\Par_X$ the set of the parabolic points for an
 action on $X.$ If $G\act X$  is a $32$-action then the set of all maximal
  parabolic subgroups is called {\it peripheral structure}   for
  the action. It  consists of  finitely many conjugacy classes
  of  maximal parabolic subgroups \cite[Main theorem, a]{Ge1}.

If  $G$ is finitely generated then an affirmative answer to the
question $\s Q2$ can be easily deduced from \cite[Map
theorem]{Ge2}   (see section 5 below). Furthermore there exists a
"universal" pullback space  in this case.    Namely every
$32$-action of a finitely generated group $G$ on a compactum $X$
admits an equivariant continuous map from the Floyd boundary
$\partial_fG$ of $G$ to $X$. The space $\partial_fG$ is universal
as it does not depend on the action on $X$ (it depends on a scalar
function $f$ rescaling the word metric of the Cayley graph and a
fixed finite set of generators of $G$).

However the same method does not work if the group is not finitely
generated. One cannot use the Cayley graph since the quotient of
the set of its edges by the group is not finite (the action is not
{\it cofinite  on edges}) the condition which is needed for the
construction of the above map. Replacing the Cayley graph by a
relative Cayley graph changes the situation since the latter graph
depends   on the $32$-action of $G$ on a compactum $X$. Indeed the
vertex set of the graph contains the parabolic points for the
action $G\act X.$ This problem turns out to be crucial since the
answer to the question $\s Q2$ is negative in general. We show in
the following theorem that a counter-example exists already in the
case of free groups of countable rank.

\medskip
 \noindent {\bf Theorem I} (Proposition \ref{freeinf}).  {\it The
 free group $F_{\infty}$ of countable rank admits
two  $32$-actions not having a   pullback space.}

\medskip

We note that this is a rare example when certain properties
 of the relatively hyperbolic groups are true for finitely generated
groups and are false for non-finitely generated (even countable)
groups.

Our next goal is to provide necessary and sufficient conditions
for two $32$-actions of a  group to have a common pullback space.

  The following theorem is the main result of the paper.

  \medskip

  \noindent {\bf Theorem II} (Theorem \ref{suffcond}, Theorem
  B).
 {\it Two $32$-actions of $G$ on
compacta $X$ and $Y$ with peripheral structures $\P$ and $\Q$
admit a pullback space $Z$ if and only if one of the following
conditions is satisfied:

\begin{itemize}

\item [1)] $\s C(X,Y):$ every element $P\in\P$ acts
$2$-cocompactly on its limit set in $Y$\rm.

\item [2)] $\s C(Y,X):$ every element $Q\in\Q$ acts
$2$-cocompactly on its limit set in $X$\rm.
\end{itemize}
}

\bigskip

Here are several remarks about the theorem. As an immediate
corollary we obtain that $\s C(X,Y)$ is equivalent to $\s C(Y,X).$
This statement seems to be new even in the finitely generated
case. It
  follows from Theorem II that a parabolic subgroup $H$ for
a $32$-action of a finitely generated group $G$ acts 2-cocompactly
on its limit set for every other such action of $G.$ Theorem I
implies that this is not true if $G$ is not finitely generated
(see Corollary \ref{equivac}.f).

 The peripheral
structure $\Rc$ for the pullback action on $Z$ is given by the
system of subgroups $\Rc=\{Q\cap P\ :\ P\in \P, Q\in\Q, \vert
P\cap Q\vert=\infty\}.$ In particular Theorem II provides a
criterion when the system $\Rc$ is a peripheral structure for some
relatively hyperbolic action of $G$ (Corollary \ref{equivac}.a).

The proof of Theorem II uses several intermediate results which
occupy first sections  of the paper and which have  independent
interest. We will now briefly describe  them.

In section 3 we study an analog of the geodesic flow introduced by
M.~Gromov in the case of hyperbolic groups. If the group $G$
admits a 32-action on a compactum $X,$ then there exists a
connected graph $\Gamma$ such that $G$ acts properly and
cofinitely on the set of edges $\G^1$ of $\G$ \cite[Theorem
A]{GePo2}. The set of vertices $\G^0$ of $\G$ is $\Par_X\sqcup G.$
The union $\ti X=X\cup \G^0=X\sqcup G$ admits a Hausdorff topology
whose restriction on $X$ and on $G$ coincide with the initial
topology and the discrete topology respectively, and $G$ acts on
$\ti X$   3-discontinuously  \cite[Proposition 8.3.1]{Ge2}. The
action is also 2-cocompact (Lemma \ref{quotient}). We call the
space $\ti X$ {\it attractor sum}  of $X$ and $G$.

Consider the space of  maps $\ga:\Z\to\ti X$ for which there exist
$m,n\in \Z\cup\{\pm\infty\}$ such that $\ga$ is constant on one or
both (possibly empty) sets
 $]-\infty, m], [n, +\infty[$ and is geodesic in
$\G^0$ outside of these sets. We call such a map  {\it eventual
geodesic} and denote by $\eG$ the space of all eventual geodesics.
We prove in section 3 (Proposition \ref{closes}) that $\eG$ is
closed in the space of maps $X^\Z$ equipped with the Tikhonov
topology. Then we show that the boundary map $\partial :\eG\to \ti
X^2$ is continuous at every non-constant eventual geodesic
(Proposition \ref{conv}). In particular we show that every two
distinct points of $\ti X$ can be joined by a geodesic (Theorem
\ref{exgeod}). This allow us to consider the convex hull ${\rm
Hull}(B)$ of a subset $B\subset \ti X$ which the union of the
images of all  geodesics in $\ti X$  with the endpoints in $B.$ We
prove that ${\rm Hull}(B)$ is closed if $B$ is.

We extensively use   so called {\it visibility property} of the
uniformity of the topology of $\ti X,$ that is for every two
disjoint closed subsets $A$ and $B$ of $\ti X$ there exists a
finite set $F\subset\G^1$ such that every  geodesic with one
endpoint in $A$ and the other in $B$ contains an edge in $F.$ At
the end of the section using the group action we show that the
space $\ti X$ cannot contain geodesic horocycles, i.e. non-trivial
 geodesics whose endpoints coincide.

 In Section 4 we study properties of subgroups of a group  acting
 $3$-discontinuously   on a compactum $X$.
 According to Bowditch \cite{Bo2} a subgroup $H$ of a group $G$   is called
\it dynamically quasiconvex \rm if for every neighborhood $\b u$
of the diagonal $\bold\Delta X$ of $X^2 = X{\times}X$ the set
$\{g{\in}G:$   $(g\bold\Lambda_XH)^2\not\subset\b u\}/H$ is
finite.

\vskip3pt Using the results of Section 3 we obtain here  the
following theorem.

\medskip

\noindent {\bf Theorem III} (Theorem \ref{infquas}). For a
$32$-action $G{\curvearrowright}X$ a subgroup $H < G$ is
dynamically quasiconvex if and only if its action on
$\bold\Lambda_XH$ is $2$-compact\rm.

\medskip

The proofs of  Theorems I and II are given in Section 5. They
 use the results of the previous sections.

  In the last section we
 provide a list of
corollaries of the main results (Corollary \ref{equivac}).

\bigskip {\bf Acknowledgements.} During the work on the paper both
authors were partially supported by the ANR grant ${\rm
BLAN}~2011\ BS01\ 013\ 04$ ''Facettes des groupes discrets''. The
first author is also thankful to   the CNRS for providing  a
research fellowship for his stay in France.

\section{Preliminaries}

\subsection{Entourages and Cauchy-Samuel completions}  We  recall some well-known notions
 from the
general topology. For further references see \cite{Ke}.

Let $X$ be a set. We denote by $\bsn$  the quotient
 of the   product space
$\underbrace{X{\times}\dots{\times}X}_{\mbox{$n$ times}}$ by the
action of the permutation group on $n$ symbols. We regard the
elements of $\bsn$ as non-ordered $n$-tuples. Let  $\Th^nX$ be the
  subset of $\bsn$ whose elements are non-ordered $n$-tuples
  whose components are all distinct. Denote
$\De^nX=\bsn\setminus \Th^nX.$

 An {\it entourage} is a
neighborhood  of the diagonal $\bde2=\{(x,x)\ :\ x\in X\}$ in
$\bs2.$ The set of  entourages of $X$ is denoted by $\ent X.$ We
use the bold font to denote entourages. For $\bu\in \ent X$ a pair
of points $(x,y)\in X^2$ is called $\bu$-small  if $(x,y)\in \bu.$
Similarly a set $A\subset X$ is $\bu$-small if $ {{\mathsf S}^2A}
\subset\bu.$ Denote by ${\rm Small }(\bu)$ the set of all
$\bu$-small subsets of $X.$

For an entourage $\bu$ we define its power $\bu^n$ as follows:
$(x,y)\in\bu^n$ if there exist  $x_i\in X$ such that $(x_{i-1},
x_i)\in \bu\ (x_0=x, x_n=y, i=1,...,n-1).$ We denote by
$\sqrt[n]\bu$ an entourage $\bv$ such that $\bv^n\subset\bu.$

A filter $\U$ on $\bs2$ whose elements are entourages is called
{\it uniformity} if
$$\foral\bu\in\U\ \exist \bv\in\U : \bv^2\subset\bu.$$

A uniformity $\U$ defines the $\U$-topology on $X$ in which every
  neighborhood  of a point has a $\bu$-small subset containing the
  point for some $\bu\in\U$. A pair $(X,\U)$ of a set $X$ equipped
  with an uniformity $\U$ is
called {\it uniform} space. A {\it Cauchy filter} $\F$ on the
uniform space $(X,\U)$ is a filter such that $\foral \bu\in \U\ :
\F\cap\s{Small}(\bu)\not=\emptyset.$ A space $X$ is {\it complete}
if every Cauchy filter on $X$ contains all neighborhoods of a
point. The uniform space $(X, \U)$ admits a completion $(\overline
X,\overline\U)$ called {\it Cauchy-Samuel} completion whose
construction is the following. Every point of $\overline X$ is the
minimal Cauchy filter $\xi$. For every $\bu\in \U$ we define an
entourage $\overline\bu$ on $\overline X$ as follows:
$$\overline\bu=\{(\xi,\eta)\in {\mathsf S}^2\overline X : \xi
\cap\eta\cap \s{Small}(\bu)\not=\emptyset\}.\hfill\eqno(1)$$ The
uniformity $\overline\U$ of $\overline X$ is the filter generated
by the entourages $\{\overline \bu : \bu\in\U\}.$ We note that the
completion $(\overline X, \overline \U)$ is {\it exact}  \cite
[II.3, Th\'eor\`eme 3]{{Bourb}}: $$\foral a, b\in\overline X \
a\not=b\ \exist \overline\bu\in\overline\U :
(a,b)\not\in\overline\bu.$$

If $X$ is a compactum then the filter of the neighborhoods of the
diagonal $\De^2X$ is the unique exact uniformity $\U$ consistent
with the topology of $X$, and $X$ equipped with $\U$ is a complete
uniform space \cite[II.4, Th\'eor\`eme 1]{Bourb}.

\subsection{Properties of ($32$)-actions of groups}
Let $X$ be a compactum, i.e a compact Hausdorff space, and $G$ be
a group acting 3-discontinuously on  $X$ (convergence action).
Recall that the limit set, denoted by $\La_XG$ (or $\La G$ if $X$
is fixed), is  the set of accumulation (limit) points of any
$G$-orbit in $X.$

 The action $G$ on $X$ is said to be {\it minimal} if $X=\La G.$

The action $G\act X$ is   {\it elementary }  if $\vert\La
G\vert<3$. If the action is not elementary then   $\La G$ is a
perfect set \cite{Tu2}. If $G$ is non-elementary then $\La G$
 is the minimal non-empty closed subset of $X$ invariant under $G$.

An elementary action of a group $G$ on $X$ is called  {\it
  parabolic} (or {\it trivial})
if $\La_XG$ is a single point. A  point $p\in \La_XG$ is {\it
parabolic} if its   stabilizer ${\rm St}_Gp$ is a maximal
parabolic subgroup fixing $p$. The set of parabolic points for the
action on $X$ is denoted by  $\Par_X.$

A parabolic fixed point $p\in \La G$ is called {\it bounded
parabolic} if the quotient space\hfil\penalty-10000 $(\La
G\setminus \{p\})/{\rm St}_G p$ is compact.

We will use an equivalent reformulation of the convergence
property in terms of {\it crosses}. A cross $(r, a)^\times\in
X\times X$ is the set $r{\times}X\cup X{\times} a$ where $(r,a)\in
X\times X.$ By identifying every $g\in G$ with its graph one can
show that $G$ acts 3-discontinuously on $X$ if and only if all the
limit points of the closure  of $G$ in $ X\times X$ are crosses
\cite[Proposition P]{Ge1}. The points $a$ and $r$ are called
respectively {\it attractive} and {\it repelling} points (or
attractor and repeller).

 A point $x\in\La G$ is {\it conical} if there is an infinite
 set $S \subset G$ such that
for every $y\in X\setminus \{x\}$ the closure of the set
$\{(s(x),s(y)) :  s \in S\}$  in $X^2$ does not intersect the
diagonal $\De^2 X.$

  A group $G$
acting on the space $X$  acts  on the set of entourages $\ent X.$
 For $\bu\in\ent X$ we denote by $g\bu$ the set
$\{(x,y)\in X^2 : g^{-1}(x,y)\in\bu\}$ and by $G\bu$ the $G$-orbit
of $\bu.$ We will say that the orbit $G\bu$ is {\it generating} if
it generates $\ent X$ as filter.

An action $G\act X$ is {\it 2-cocompact} if $\Theta^2X/G$ is
compact. Suppose that a group $G$ admits a
   3-discontinuous  and 2-cocompact non-parabolic minimal action
($32$-action)   on a compactum $X$.  Then every point of $X$ is
either   a bounded parabolic   or  conical \cite[Main
Theorem]{Ge1}. P.~Tukia showed that if   $X$ is metrisable then
the converse statement is  true \cite[Theorem 1C, (b)]{Tu2}.

Let $\G$ be  a graph. We  denote by $\G^0$ and $\G^1$  the set of
vertices and edges of $\G$ respectively. Recall that an action of
$G$ on $\Gamma$ is {\it proper on edges} if the stabilizer ${\rm
St}_{\G}e$ of every edge $e$ in $\G$ is finite. The action $G\act
\G$ is called {\it cofinite} if $\vert\Gamma^1/G\vert < \infty$.

 According to
B.~Bowditch \cite{Bo1} a graph $\Gamma$ is called \it fine \rm if
for any two vertices the set of simple arcs of fixed length
joining them is finite.

It is shown in \cite[Theorem 3.1]{GePo2} that if $G$ admits a
non-parabolic $32$-action on a compactum $X$ then Bowditch's
condition of relatively hyperbolicity is satisfied. This means
that there exists a connected fine and hyperbolic graph $\G$ acted
upon by $G$ cofinitely and properly on edges. Every vertex of $\G$
is either an element of $G$, or belongs to the set of parabolic
points ${\mathsf Par}_X$.

 Consider    the union of two topological spaces $\ti X=X\sqcup G=X\cup \G^0$
 where $G$ is equipped with the discrete topology.  By \cite[Proposition 8.3.1]{Ge2}
 $\ti X$
 admits a unique compact Hausdorff
topology  whose restrictions on $X$ and $G$ coincide with the
original topologies of $X$  and $G$, and the action on $\ti X$ is
3-discontinuous (the description of this topology see in
Proposition \ref{pullback})

Following \cite{Ge2} we call the space$\ti X$ {\it attractor sum}
of $X$ and  $\G$.

The action on $\ti X$ is also 2-cocompact. Indeed by assumption
the action of $G\act \Theta^2 X$ is cocompact. So there exists  a
compact fundamental set $K\subset \Theta^2 X$. Hence $\ti K=K\cup
(\{1\} \times (\ti X\setminus\{1\}))$ is a compact fundamental set
for the action on $\Theta^2\ti X.$  Therefore the action $G\act
\ti X$ is a $32$-action. We summarize all these facts in  the
following lemma.

\begin{lem}\label{cocomext}  \cite{Ge2},  \cite{GePo2}.  Let $G$ admits a   non-parabolic  $32$-action
on a compactum X. Then there exists  a connected, fine and
hyperbolic graph   $\G$ acted upon by $G$
  properly  and cofinitely on edges.
Furthermore  $G$ acts 3-discontinuously and 2-cocompactly
  on the attractor sum $\ti X=X\cup\G^0$ and
   $\G^0=G\cup {\mathsf Par}_X$ is the set of all non-conical points
  for the action on $\ti X$. \bx
\end{lem}

 We will consider the entourages
 $\bu\in \ent\ti X$ on the attractor sum $\ti X$ as well as
 their restrictions on $\G$ and on $X.$

\medskip

Following Bowditch \cite{Bo1} for a fixed group $G$  a
$G$-invariant set $M$ is called {\it connected $G$-set} if  there
exists a connected graph $\G$ such that $M=\G^0$ and the action
$G\act\G^1$ on edges is proper and cofinite.

Recall some more definitions. An entourage $\bu$ on a connected
$G$-set $M$ is called   {\it perspective} if for any pair
$(a,b)\in M\times M$ the set $\{g\in G\ \:\ g(a,b)\not\in\bu\}$ is
finite.

An entourage $\bu$ given on a connected $G$-set $M$ is called {\it
divider} if there exists a finite set $F\subset G$ such that
$({\bu} _F)^2\subset \bu$ where $\bu_F=\cap_{f\in F}f\bu.$

We say that uniformity $\U$ of a compactum $\ti X$ is {\it
generated by an entourage $\bu$} if it is generated as a filter by
the orbit $G\bu $.

\begin{lem}\label{perdiv}\cite[Proposition 8.4.1]{Ge2}.
If a group $G$ acts $3$-discontinuously and $2$-cocompactly on a
compact space $X$ then the uniformity $\U$ on the compactum $\ti
X=X\cap G$ is generated by a perspective divider $\bu.$
\end{lem}

The following result describes the opposite  way  which starts
from a perspective
  divider  on a connected $G$-set $M=\G^0$  and gives   a
  $32$-action on the compactum $\ti X=X\cup \G$ where $X$ is a
  "boundary" of $\G$.

  \begin{dfn}\label{visprop} Let $e\in\G^1$ be an edge.
A pair of vertices $(a,b)$ of $\G$ is called $\bu_e$-small  if
there
  exists  a geodesic in $\G$ with endpoints $a$ and $b$
which does not contain    $e.$

A uniformity $\U^0$ on $M=\G^0$ has a {\it visibility property} if
for every entourage $\bu^0\in \U^0$ there exists a finite set of
edges $F\subset \Gamma^1$ such that $\bu_F=\cap\{\bu_e\ \vert\
e\in F\}\subset\bu^0.$
\end{dfn}

\noindent The   following  lemma describes the completion $\ti X$
mentioned above and will be often used in the paper.

\begin{lem}\label{compldense} \cite[Propositions 3.5.1, 4.2.2]{Ge2}. Suppose that a group $G$ acts
on a connected  graph $\G$
  properly  and cofinitely on edges. Let $\W^0$ be a uniformity
  on
  $\G^0$  generated by a
  perspective
  divider. Then  $\W^0$ has the visibility property. Furthermore
  the Cauchy-Samuel completion $(Z,\W)$ of the uniform space $(\G^0, \W^0)$ admits
   a $32$-action
  of $G$.
\end{lem}

\medskip

Let $\G$ be a connected graph. We now recall the definition of the
{\it Floyd  completion (boundary)}  of $\G$ mentioned in the
Introduction (see also \cite{F}, \cite {Ka}, \cite{Ge2}).

A function $f:\Bbb N\to\Bbb R$ is said to be a (Floyd) \it scaling
function \rm if $\sum_{n\geqslant0}f_n<\infty$ and there exists a
positive $\lambda$ such that $1\geqslant
f_{n+1}/f_n\geqslant\lambda$ for all $n{\in}\Bbb N$.

Let $f$ be a scaling function and let $\Gamma$ be a connected
graph. For each vertex $v{\in}\Gamma^0$ we define on $\Gamma^0$ a
path metric $\d_{v,f}$ for which the length of every edge $e\in
\Gamma^1$ is $f(d(v, e))$. We say that $\d_{v,f}$ is the \it Floyd
metric \rm(with respect to the scaling function $f$) \it based \rm
at $v$.

When $f$ and $v$ are  fixed we write $\d$ instead of $\d_{v,f}$.

One verifies that $\d_u/\d_v\geqslant\lambda^{\s d(u,v)}$ for
$u,v{\in}\Gamma^0$. Thus the Cauchy completion $\overline\Gamma_f$
of $\Gamma^0$ with respect to $\d_{v,f}$ does not depend on $v$.
The \it Floyd boundary \rm is the space $\partial_f\Gamma
=\overline\Gamma_f\setminus \Gamma^0$. Every $\s d$-isometry of
$\Gamma$ extends to a homeomorphism
$\overline\Gamma_f\to\overline\Gamma_f$. The Floyd metrics extend
continuously onto the Floyd completion $\overline\Gamma_f$.

In the particular case when $\G$ is a Cayley graph of $G$ we
denote by $\partial_fG$ its Floyd boundary  or by $\partial G$ if
$f$ is fixed.

\section{Geodesic flows on graphs} In this section we study the
properties of geodesics on a class of  graphs. Let $\G$ be a
connected graph. We will assume that $\G^0\subset \ti X$ for a
compactum $\ti X.$ Let $\U$ be the  uniformity consistent with the
topology of $\ti X.$ Since $\ti X$ is Hausdorff the uniformity
$\U$ is exact. In this section we will always admit the following.

\medskip

\noindent {\bf Assumption.}
 The uniformity $\U$ has the visibility property on $\G^0.$

\medskip

The most of the material of the section does not relate to  any
group action. However the only known example when the above
assumption is satisfied is the case when a compactum $X$ admits a
$32$-action of a group $G$ and $\ti X=X\cup \G$ is the attractor
sum   (see Lemma \ref{compldense}).

A path in $\G$ is a map $\ga : \Z\to \G$ such that $\ga\{n, n+1\}$
is either an edge of $\G$ or a point $\ga(n)=\ga(n+1).$ A path
$\ga$ can contain a "stop" subpath, i.e. a   subset $J$ of
consecutive integers   such that $\ga\vert_J\equiv{\rm const}.$

For a finite subset $I\subset\Z$ of consecutive integers we define
the boundary $\partial(\gamma\vert_I)$ to be $\gamma(\partial I)$.
We extend naturally the meaning of $\partial\gamma$ over the
half-infinite and bi-infinite paths in $\G^0\subset\ti X$ in the
case if the corresponding half-infinite branches of $\gamma$
converge to points in $\ti X$. The latter one means that for every
entourage $\bv\in\U$ the set $\gamma\vert_{[n,\infty[}$ is
$\bv$-small for a sufficiently big   $n$.

\medskip

\begin{lem}\label{conv}
Every half-infinite geodesic ray  $\ga:[0, \infty[\to\G$ converges
to a point in $\ti X.$
\end{lem}

\proof  Fix an entourage $\bv\in\U.$  By the visibility property
  there exists a finite set of edges $F\subset\Gamma^1$
such that $\displaystyle\bu_F=\bigcap_{e\in F}\bu_e\subset \bv.$
Since $\gamma$ is a geodesic, the ray $\gamma\vert_{[n_0,\infty[}$
does not contain $F$ for some  $n_0\in \N$. So
$\gamma\vert_{[n_0,\infty[}$ is $\bu_F$-small and  therefore
$\bv$-small.\bx

\begin{dfn}
\label{evgeod}

A path $\gamma:I\to \G$ is an eventual geodesic if it is either a
constant map, or each   its maximal stop-path  is  infinite and
outside of its maximal stop-paths $\ga$ is a geodesic in $\G^0.$

 The set of eventual geodesics in $\G$ is denoted
by   $\rm EG(\G).$

\end{dfn}

 The image of every eventual
geodesic in $\G$ is either a geodesic (one-ended, or two-ended or
finite) or a point.

\begin{prop}\label{closes} The space $\eG$ is closed in the space
of maps $\ti X^\Z$ equipped with the Tikhonov topology.
\end{prop}

\proof  Let $\al$ belongs to the closure $\overline{\eG}$ of $\eG$
in $\ti X^\Z$. If $\al$ is a constant then $\al\in \eG$ and there
is nothing to prove. So we assume that $\al\in\overline{\eG}$ is a
non-trivial map. The proof follows from the following three lemmas
having their own interest.

\medskip

\noindent {\bf Lemma \ref{closes}.1}\ {\it If
$\al(n)\not=\al(n+1)$ then there exists a neighborhood $O$ of
$\al$ in $\X^{\mathbb Z}$ such that $\gamma(n)=\al(n)$,
$\gamma(n+1)=\al(n+1)$ for every $\gamma\in O \cap \eG$. In
particular $\{\al(n),\al(n+1)\}$ is an edge of $\G.$}

\medskip

{\it Proof of the Lemma.} Since $\U$ is an exact uniformity there
exists $\bu\in \U$ such that $(\al(n),\al(n+1))\notin\b u\in\c U$
(see Section 2.1). Since $\U$ is a uniformity  $\exist \bv\in\U\
:\ \bv^3\in\U$ and $\bv\in\U$.   By the visibility property
 there exists a finite set $F\subset
\G^1$ such that $\bu_F\subset\bv$. Let $F^0$ denote the set of
vertices of the edges in $F$. Let $O_n$ be a $\b v$-small
neighborhood of $\al(n)$ disjoint from $F^0{\setminus}\{\al(n)\}$
and   $O_{n+1}$ be a neighborhood of $\al(n+1)$ defined in the
same way.

Set  $O=\{\gamma\in\X^{\mathbb Z}:\gamma(n)\in O_n\}$. If
$\gamma\in O$ then $(\ga(n), \ga(n+1))\not\in \bv$ and
$\ga(n)\not=\ga(n+1).$ If in addition $\ga\in \eG$ then
$\ga\vert_{[n, n+1]}$ is a geodesic and necessarily
$\{\gamma(n),\gamma(n+1)\}\in F$. By the definition of $O_n$ and
$O_{n+1}$ we obtain $\gamma(n)=\al(n)$, $\gamma(n+1)=\al(n+1)$.\bx

\medskip

 \noindent {\bf Lemma \ref{closes}.2. {\it Every maximal stop for $\al$ is
infinite.}

\medskip

\noindent \it Proof of Lemma\rm. By contradiction assume that
$J\subset\Z$ is a finite maximal stop for $\al$. By Lemma 1 this
means that, for $m,n\in\mathbb Z$ such that $n-m{\geqslant}3$ one
has $\al(m)\ne\al(m+1)$, $\al(n-1)\ne\al(n)$ and
$\al(k)=b\in\Gamma^0$ for $m<k<n$. Moreover, by Lemma 1 there
exists a neighborhood $O$ of $\al$ in $\X^{\mathbb Z}$ such that
if $\gamma \in O$ then $\gamma(m)=\al(m)$, $\gamma(m+1)=\al(m+1)$,
$\gamma(n-1)=\al(n-1)$, $\gamma(n)=\al(n)$. Since $\al$ belongs to
the closure of $\eG$ such $\gamma$ does exist. But this implies
that $\gamma(m+1)=\gamma(n-1)=b$. As $\gamma$ is an eventual
geodesic the interval $\{k\in\mathbb Z:m < k < n\}$ is a stop for
$\gamma$ so either $\al(m)=\gamma(m)=b$ or $\al(n)=\gamma(n)=b$
contradicting the maximality of $J.$\bx

\medskip

\noindent {\bf Lemma \ref{closes}.3}. {\it If $J$ is a finite
  subset of consecutive integers and $J$ does not contain
    stops for $\al$ then $\al|_J$
is geodesic.}

\medskip

 \it Proof of the Lemma\rm. By the hypothesis each two
consecutive values of $\al|_J$ are distinct. Since $J$ is finite
by Lemma 1 there exists a neighborhood $O$ of $\al$ in
$\X^{\mathbb Z}$ such that if $\gamma \in O \cap \eG$ then
$\gamma|_J=\al|_J$. Since $\al\in\overline{\eG}$ such $\gamma$
does exist.\bx

\medskip

 It follows from lemmas 2 and 3 that $\al\in\eG$.  The proposition
 is proved.

\medskip

\begin{cor}\label{clbal} For every $n\geqslant 0$ the ball $\s B_n  q$ of radius $n$
in $\Gamma^0$ centered at any $ q\in\Gamma^0$ is closed in $\X$.
\end{cor}

\medskip

 \it Proof\rm. Let $ p\in\overline{\s B_n  q} \setminus \s
B_n  q$. For a neighborhood $O$ of $ p$ in $\X$ let
$\gamma_O:\{k\in\mathbb Z:0{\leqslant}k{\leqslant}n\}\to\Gamma^0$
be a geodesic joining $  q$ with a point in $O$. We make each such
$\gamma_O$ an eventual geodesic by extending it by constants.
Since $\X^{\mathbb Z}$ is compact there is an accumulation point
$\gamma$ for the set of all $\gamma_O$. By Proposition
\ref{closes} $\gamma\in\eG$. Since the projections $\X^{\mathbb
Z}\to\X$ are continuous and $\X$ is Hausdorff we have $\gamma(0)=
q$, $\gamma(n)= p$. Since $\gamma$ is eventual geodesic we have $
p\in\Gamma^0$.\bx

\medskip

\begin{cor}\label{opnev}
For every finite path $l=\{a_1, ..., a_n\}\subset \G$ the set
$${(\eG)_l}=\{\ga\in \eG : \ga(I)=l,\ I\subset Z, \
\ga(-\infty)\not=a_1,\ \ga(\infty)\not=a_n\}$$ is open.
\end{cor}
 \proof By  the proof of Proposition \ref{closes} $\gamma$ admits a neighborhood
$O\subset\eG$ such that $\foral \la\in O\ \ga\vert_l=\la\vert_l$.
\bx

\medskip

By Lemma \ref{conv} for    a half-infinite geodesic ray
$\ga:\Z_{>0}\to X$   $\displaystyle\lim_{t\to +\infty}\ga(t)$
exists.  The following proposition refines   this statement.
\begin{prop}\label{bcont}
The boundary map $\displaystyle \partial :\eG\to \ti X^2$ where
$\displaystyle
\partial :\ga\to
\partial\ga=\{\lim_{t\to-\infty}\ga(t), \lim_{t\to+\infty}\ga(t)\}$ is continuous at every
non-constant eventual geodesic.
\end{prop}

\proof For  $\al\in \eG$ we denote by $\al_{-\infty}$ and
$\al_{+\infty}$ the limits $\displaystyle\lim_{t\to-\infty}\al(t)$
and $\displaystyle\lim_{t\to+\infty}\al(t)$ respectively. We will
prove that both coordinate functions   $\displaystyle\pi_- :
\al\to
  \al_{-\infty}$ and $ \pi_+ : \al\to\al_{+\infty}$ are
   continuous for every non-constant
geodesic $\al.$

Fix $\al\in \eG$ and suppose  that $a=\al_{+\infty}.$ We need to
prove that for every small neighborhood $U$ of $a$ there exists a
neighborhood of $N_\al\opn\eG$  such that one of the endpoints of
every eventual geodesic belonging to $N_\al$  is in $U.$

\underline{Case 1.}  $\al (n)\ (n\geq 0)$ are all distinct.

  Let $U$ be a closed
neighborhood of $a$ such that $b=\al(0)\not\in U.$ Choose a
"smaller" closed neighborhood $V$ of $a$ such that the interior
$\accentset{\circ}{U}$ of $U$ contains $V.$ By the exactness of
the uniformity $\U$ there exists an entourage $\bv\in\U$ such that
$\bv\cap (U\times V)=\emptyset.$ Then, by the visibility property,
we have $\bu_F=\cap_{e\in F} \bu_e\subset \bv$ for some finite set
$F$ of edges of $\G.$ So every eventual geodesic $\ga$ passing
from $b\in U'=\ti X\setminus U$ to $V$ contains an edge from $F.$
Denote by $d$    the diameter (in the graph distance) ${\rm diam}
F={\rm max} \{d(a_i, a_j)\ :\ a_i\in F^0\}$ of the set of vertices
$F^0$ of $F.$ Since $\G$ is connected and $F$ is finite $d$ is
finite. Since $\{\al(n)\}_n$ converges to $a$ there exists  $m$
such that $\al(m+i) \in \accentset{\circ}{V}\cap \G^0$ for $
i=0,...,d+1.$

Consider the following set:
$$N_\al=\{\ga\in\eG\ :\
\{\ga(m),...,\ga(m+d+1)\}\subset
{\accentset{\circ}{V}}\cap\G^0\}.\hfill\eqno(1)$$

 We have $N_\al\not=\emptyset$ as $\al\in N_\al.$
 Furthermore  $N_\al$ is open in $\eG.$  Indeed ${\accentset{\circ}{V}}\cap\G^0$
  is open in $\G^0$ so
the condition $\ga(t)\in {\accentset{\circ}{V}}\cap\G^0$ defines
an open subset of $\eG$ and $N_\al$ is the intersection of
finitely many such subsets.

Let $\ga\in N_\al.$ We claim that $\ga$
  cannot quit   ${U}\cap \G$. Indeed if not
   then $\ga$ contains a finite subpath $\ga'$ which passes from $U$ to
   $V$, then it   passes through at least $d+1$
    distinct consecutive vertices $\ga(i)\in  {\accentset{\circ}{V}}\ (i=m,...,m+d+1),$
      and after  it goes back to $U'=\ti X\setminus U.$
Assuming    $\ga'$ to be the minimal subpath of $\ga$ having these
properties, by Lemma \ref{closes}.2 we obtain that $\ga'$
  a geodesic of length at least $d+1.$
Then there is a couple   $(i,j)$ of indices such that $i \leq m,\
j\geq m+d+1,\
  \{q_i=\ga'(i),q_j=\ga'(j)\}\subset F^0$ and $d(q_i, q_j) \geq d+1 > {\rm diam} (F)$ what is impossible.

\medskip

\underline{Case 2.}   $a\in \G^0$.

A similar  argument works in this case too. We can assume that
$a=\al(0)$. Up to re-parametrisation of $\al$ we can assume that
$\foral t\geq 0\ :\ \al(t)=a$ and $\al(-1)=b\not=a.$   Consider
two disjoint closed neighborhoods $U$ and $V$ of $a$ such that
$V\subset {\accentset{\circ}{U}}$ and $b\not\in U$. As above there
exists a finite set $F\subset\G^1$ such that every eventual
geodesic passing from $V$ to $U$ (or vice versa) contains an
  edge of $F$. Let $d={\rm diam} F.$

 For $k\geq d+1$  put

$$N_\al= \{\ga\in \eG\ :\ \ga(-2)=\al(-2), \ga(-1)=\al(-1), \{\ga(0), ..., \ga(k)\}\subset{\accentset{\circ}{V}}
\}.\hfill\eqno(1')$$

The set $N_\al$ is non-empty as $\al\in N_\al$. We have
$N_\al=N_1\cap N_2$ where $N_1=\{\ga\in\eG\ :\ \ga(-2)=\al(-2),
\ga(-1)=\al(-1)\}$ and $N_2=\{\ga\in\eG\ :\{\ga(0), ...,
\ga(k)\}\subset{\accentset{\circ}{V}} \}.$ Since $\al(-2)\not
=\al(-1)$   by Corollary \ref{opnev} $N_1$ is open. The set $N_2$
is also open so is $N_\al.$

Every $\ga\in N_\al$ passes from the point $b$ to $V$ and admits a
geodesic sub-interval in $V$ of length at least $d+1.$ So by the
argument of   Case 1 we have that  $\ga_{+\infty}\in U.$

 We have proved that the
coordinate functions $\pi_-$ and $\pi_+$ are continuous. Therefore
$\partial=(\pi_-,\pi_+)$ is continuous at every non-trivial
eventual geodesic.\bx

\medskip

\noindent {\bf Remark.}  Note that the map $\partial$ is not
continuous at any constant eventual geodesic $\al(t)=a\in  \G^0
(t\in \Z).$ Indeed the eventual geodesics of the form
$\displaystyle\{\ga\in\eG\ :\  \ga\vert_{]-\infty,n]}\equiv a, \
\ga_{+\infty}=b\not=a\}$ converge to $\al$ when $n\to \infty$ in
the Tikhonov topology but $\pi_+(\ga)\not=a.$

\medskip

\begin{theor}\label{exgeod}
For every two distinct  points $p,q\in \ti X$ there exists an
eventual geodesic $\gamma$ such that $\gamma_{-\infty}{=} p$,
$\gamma_{+\infty} =  q$.
\end{theor}

\proof  If both $ p$ and $ q$ belong to $\Gamma^0$ then the
assertion follows from the connectedness of $\Gamma$. So we
suppose that at least one of the points does not belong to
$\Gamma^0$. We first consider the case when  both $ p$ and $ q$ do
not belong to $\Gamma^0$. Then we explain how to modify the
argument  when exactly one of the points belongs to $\Gamma^0$.

Let $P_0,Q_0$ be two closed disjoint neighborhoods of ${p,q}$
respectively. By the exactness of $\U$ there exists an entourage
$\b u_0\in\U$ such that $\bu_0\cap P_0{\times}Q_0=\emptyset$. By
the visibility property there exists a finite set
$F{\subset}\Gamma^1$, such that $\b u_F\subset \b u_0$. Let $\{(
a_i, b_i):i{=}0,1,\dots,m\}$ be the list of \bf ordered \rm pairs
such that $\{ a_i, b_i\}{\in}F$.

For closed neighborhoods $P,Q$ of ${p,q}$  contained in $P_0,Q_0$
respectively let

 $$W_{i,P,Q} =\{\gamma\in\eG:\gamma_{-\infty}\in
P,\gamma_\infty\in Q,\gamma(0)= a_i,\gamma(1) =
b_i\}.\hfill\eqno(2)$$

We claim that $W_{i,P,Q}$ is closed. Indeed suppose $\gamma\notin
W_{i,P,Q}$. If $\gamma$ is not a constant then Proposition
\ref{closes} implies that the opposite of every condition in (2)
defines an open subset in $\eG.$ Then their finite intersection is
open. If $\gamma\equiv c$ is a constant then either
$\ga(0)\not=a_i$ or $\ga(1)\not= b_i$. So there is an open
neighborhood $N_\ga\subset\eG$ such that every $\beta\in N_\ga$
satisfies $\{\beta\in\eG:\beta(0)\ne a_i\}$ or
$\{\beta\in\eG:\beta(1)\ne b_i\}$ respectively. In each case
$N_\ga\cap W_{i,P,Q}=\emptyset$ and the claim follows.

There exists $i\in\{0,...,m\}$ such that all $W_{i,P,Q}$ are
nonempty. Indeed if not then for every $i\in\{1,...,m\}$ there
exist neighborhoods  $P_i$ and  $Q_i$ of $p$ and $q$ such that
$W_{i,P_i,Q_i}=\emptyset$. Then, there would no geodesic between
the closed non-empty subsets  ${\bigcap}_{i=0}^mP_i$ and
${\bigcap}_{i=0}^mQ_i$ which is impossible. Then for some $i$, say
$i=0$, the family $W_{0,P,Q}\not=\emptyset$   for all $P$ and $Q.$
So it is   a centered family of non-empty closed sets. Since $\eG$
is compact $\exist\ga\in\eG$ such that
$\ga\in\bigcap_{P,Q}W_{0,P,Q}.$ By definition of $W_{0,P,Q}$ the
point $\gamma_{+\infty}$ belongs to any neighborhood of $ q$.
Hence $\gamma_{+\infty} =  q$, and similarly $\gamma_{-\infty} =
p$.

The assertion is proved for the case $ p\notin \Gamma^0$, $
q\notin\Gamma^0$. If one of the vertices, say $ p$, belongs to
$\Gamma^0$ then we modify the definition of $W_{i,P,Q}$ as
follows: $W_{i,P,Q} = \{\gamma\in \eG : \gamma_{-\infty} =
\gamma(-n_i) =  p,\gamma_{+\infty} \in Q,\gamma(0) = a_i,\gamma(1)
= b_i\}$ where $n_i$ is the distance between $ p$ and $ a_i$. Then
the above argument works without any change.
 \bx

\medskip

\noindent Let $B\subset \ti X$ be a closed set. Define its
(eventual) geodesic hull as follows:

$$ {\rm Hull}(B)=\cup\{{\rm Im}\g  : \partial\g\subset
  B\}.\hfill\eqno(3)$$

\noindent The following lemma and its corollary will be used
further.

\begin{lem} \label{closgeod} The set

$$\C=\{\ga\in\eG\ :\ \partial\ga\in B^2\}\hfill\eqno(3)$$

\noindent is closed in $\eG.$
\end{lem}

\begin{cor} \label{hull} If $ B\subset\ti X$ is closed then
${\rm Hull}(B)$ is a closed subset.

\end{cor}

\noindent {\it Proof of the corollary.} The projection $\pi:\ga\in
\C\to \ga(0)$ is continuous by the definition of the Tikhonov
topology. Since $X^\Z$ is compact and if $\C$ is closed then
$\pi(\C)$ is closed.\bx

\medskip

\noindent {\it Proof of the lemma.}  Let us  show that
$\eG\setminus \C$ is open. Let $\al\in\eG\setminus \C$ does not
belong to $\C$.

Suppose first that $\al$ is not a constant. Then
$\partial\al\not\in B^2$. Since $B^2$ is closed in $\ti X^2$ there
exists an open neighborhood $W$ of $\partial\al$ such that $W\cap
B^2=\emptyset.$ By Proposition \ref{bcont}
$N_\al=\partial^{-1}(W)$ is  an open subset of $\eG$, so we are
done in this case.

Suppose now that $\al$ is a constant: $\al\equiv a\in\ti
X\setminus {\rm Hull}(B)$. Choose  a closed   neighborhood $U$ of
$a$ disjoint from   $B$.  By the visibility property there exists
a finite set $F\subset\G^1$ of a finite diameter $d$ such that
every eventual geodesic passing from $U$ to $B$ passes through
$F.$

Suppose by contradiction that every neighborhood $N_\al$ of $a$ in
$\eG$ intersects $\C$. Note that every such $N_\al$ contains a
non-trivial eventual geodesic $\ga\in \eG$  since otherwise $U\cap
B\not=\emptyset$ which is not possible. So  for any $k\in \Z_+$ we
can find a non-trivial $\ga\in \C$ such that $\partial\ga\in B^2$
and $\ga(i)\in U\ (i=0,1,2...,k)$. Then $\ga$  contains an edge
from $F$ on its way  from $B$ to $U$  and it meets $F$ again on
its way back from $U$ to $B$. So there is a geodesic subsegment
$\ga'$
 of $\ga$ passing through $k$ consecutive points
in $U$ and having its endpoints in $F.$ Therefore  for $k\geq d$
it would imply ${\rm diam} F \geq k \geq {\rm length}(\ga') > d$
which is a contradiction.\bx

\medskip

At the end of this section we obtain a description of the
endpoints of one-ended and two-ended eventual geodesics. It is the
first (and the last) time in the section when we use a group
action

 Let a group $G$ admits a non-parabolic 32-action on a compactum $X$. Then by
 Lemma \ref{cocomext} there exists a connected and fine graph $\G$ such that
 the  action on $\G$ is proper and cofinite on the edges.
We also suppose that $\Gamma^0 = G\cup\{$the parabolic points
${\in}\widetilde X\}=\{$the non-conical points of $\widetilde
X\}$. The existence of a $G$-finite $G$-set $\Gamma^1$ making
$\Gamma^0$ into the vertex set of a connected graph is proved in
\cite[Theorem 3.1]{GePo2}. We note that it is also shown in
\cite{GePo2} that the graph $\G$ is hyperbolic but we will not use
it in this section.

\begin{prop}\label{horabs}
 For every non-trivial eventual geodesic $\gamma$ one has
$\gamma_{-\infty}\ne\gamma_{+\infty}$.

\end{prop}

\noindent {\bf Remark.} In other words the proposition claims that
the  graph $\G$ does not have non-trivial eventual geodesic
horocycles.

\medskip

 \proof Since the action $G\act X$ is 3-discontinuous and
2-cocompact every limit point is either conical or bounded
parabolic \cite[Main theorem]{Ge1}. We can assume that the action
on $X$ is minimal so $X=\La_XG.$

Suppose first that $ p=\gamma_{-\infty}= \gamma_{+\infty}$ is
conical. Then there exist closed disjoint sets $A,B \subset   X$
and infinite subset $S$ of $G$ such that    for every closed set
$C\subset \ti X\setminus \{p\}$ and every closed neighborhood $\ti
B$ of $B$ in $\ti X$ disjoint from $A$ there exists a subset
$S'\subset S$ such that $\vert S\setminus S'\vert\ <\infty $,
$S'(C)\subset \ti B$ and $S'(p)\in A.$
 For an
arbitrary finite set $J \subset \mathbb Z$   we have $s\gamma(J)
\subset \widetilde B$ for all but finitely many $s \in S$. By
Lemma \ref{closgeod} the set $\partial^{-1}(\widetilde B)$ is
closed. So every accumulation point of the set $S\gamma$ belongs
to $\partial^{-1}\widetilde B^2$.   On the other hand every
eventual geodesic $s\gamma$ belongs to the closed set
$\partial^{-1}A^2$ as $\partial(s\ga)=s(p)\in A.$ So
$s\gamma\in\partial^{-1}\widetilde B^2\cap
\partial^{-1}A^2$ which is impossible as $A^2\cap \ti
B^2=\emptyset.$

Let $p = \gamma_{-\infty}{=}\gamma_{+\infty}\in\Gamma^0$. In this
case $\gamma$ can not be finite so $  p$ is a limit point for the
action $G\act\X$. Hence $p$ is bounded parabolic. Let $K$ be a
compact fundamental set for the action $\s{St}_G
p{\curvearrowright}\X{\setminus}\{ p\}$.  Let $\bu$ be an
entourage such that $\bu\cap (\{ p\}\times K) =\emptyset.$ By the
visibility property there exists a finite set $F\subset\G^1$ such
that every eventual geodesic from $p$ to $K$ contains an edge in
$F.$ For every $n\in\mathbb Z$ there exists $h\in\s{St}_G p$ such
that $h\gamma(n)\in K$. Since $\partial(h\gamma)= ({p,p}),$ the
geodesic $h\gamma$ has edges in $F$ in both sides of $h\gamma(n)$.
Hence $\s{dist}_\Gamma( p,\gamma(n))=\s{dist}_\Gamma(
p,h\gamma(n))\leqslant\s{dist}_\Gamma(p,F^0){+}\s{diam}_\Gamma
F^0<\infty\ (n\in\Z).$ So $\gamma(\mathbb Z)$ is bounded and can
not be an infinite geodesic.\bx

\medskip

\begin{prop}\label{geodray}  Let
$\gamma:[0,+\infty)\to\Gamma^0$ be an  infinite geodesic ray then
  $\gamma_{+\infty}\notin\Gamma^0$ and  it is a conical point.
\end{prop}

 \proof If the assertion were false the point $  p=\gamma_{+\infty}$ is
 bounded parabolic.
Let $K$ be a compact fundamental set for the action $\s{St}_G
p\act\X{\setminus}\{ p\}$. Since the geodesic ray is infinite the
graph distance $d(\ga(0), \ga(n)$ tend to infinity. Let $h_n$ be
an element of $\s{St}_G\g p$ such that maps $h_n(\gamma(n))\in K
(n\in\Z_+)$. By applying the ``flow map'' $\ga(n)\to \ga(n-1)$ to
the geodesic $h_n\gamma$ we obtain a  geodesic  $\gamma_n$ such
that $\gamma_n(0)\in K$. The sequence of their other endpoints
$\{h_n(\gamma(-n))$ converges to $ p$ as $n\to+\infty.$ Each
accumulation point $\al$ of the set $\{\gamma_n\ :\ n \geqslant
0\}$ in $\eG$ is an eventual geodesic whose both endpoints are
equal to $p$. Furthermore it cannot be a constant since it has a
value in $K$. We have a non-trivial $\al\in\eG$ such that
$\al_{-\infty}=\al_{+\infty}$ which contradicts Proposition
\ref{horabs}.\bx

\medskip

In the following Proposition we show that the absence of the
horocycles is equivalent to the hyperbolicity of the graph $\G$.


\begin{prop}\label{uniGen}
The uniformity $\c U^0=\c U\vert_{\G^0}$  is generated by the
collection $\{\mathbf u_e:e{\in}\Gamma^1\}$.
\end{prop}

\noindent\it Proof\rm. Let us prove that $\mathbf u_e$ is an
entourage. By Proposition \ref{bcont} the boundary map $\partial$
is continuous on the closed set
$K_{0,e}=\{\gamma\in\mathrm{EG}(\Gamma):e=\{\gamma(0),
\gamma(1)\}\}$. Hence $\partial K_{0,e}$ is closed. Its complement
in $\widetilde X^2$ is exactly $\mathbf u_e$. By Proposition
\ref{horabs} the open set $\mathbf u_e$ contains all diagonal
pairs $(p,p)$, $p{\in}\widetilde X,$ and so is an entourage.

By the visibility property for every $\bv^0\in \c U^0$ there
exists a finite set $F$ of edges for which $\bu_F=\cap_{e\in
F}\bu_e\subset \bv^0$. So  $\c U^0$ is generated by the set
$\{\bu_e\ :e\in\G^1\}$ as a filter. \hfill$\square$

\medskip
\bf Remark\rm. By [Ge2, subection 5.1] this corollary implies that
$\Gamma$ is hyperbolic. Namely each side of every geodesic
triangle $\triangle$ is contained in the metric
$\delta$-neighborhood of the union of the other sides. The
constant $\delta$ is determined as follows. Let $E_\#$ be a finite
set of edges intersecting each $G$-orbit in $\Gamma^1$. For every
$e\in E_\#$, since $\mathbf u_e$ is an entourage, there exists a
finite set $F\subset\Gamma^1$ such that $\mathbf
u_{F}^2\subset\mathbf u_e$ (this property is called \it
alt-hyperbolicity \rm in [Ge2]). It follows directly from the
definition of $\mathbf u_e$ that one can choose
$\delta=1{+}\s{max}\{\s{dist}_\Gamma(e,F(e)^0):e{\in}E_\#\}$.

This gives an independent proof of Yaman theorem without
metrisability and cardinality restrictions. Note that it uses
[GePo2] where a connected graph $\Gamma$ was constructed. It also
uses [Ge2] where the visibility property was proved.


\section {Dynamical quasiconvexity and 2-cocompactness condition.}\label{secqconvex}

\subsection{The statement of the result} Let $X$ be a compactum.
We  first restate the definition of the dynamical quasiconvexity
in terms of entourages \cite{GePo3}.

\medskip

\begin{dfn}\label{dynquas}    Let $G$ be a
discrete   group acting 3-discontinuously on a compactum $X$. A
subgroup $H$ of $G$  is said to be   dynamically quasiconvex   if
for every entourage $\mathbf u$ of $X$ the set $G_{\bold u}=\{g
G:g(\La H)\notin\s{Small}(\bold u)\}/H$ is finite.

\end{dfn}

 The aim of this Section is the following theorem giving a characterization
 of the dynamical quasiconvexity.

\vspace*{0.5cm}

\begin{thm} \label{infquas}.   Let $G$ be a group which admits
3-discontinuous and 2-cocompact non-trivial action on a compactum
$X$. Let $H$ be a subgroup of $G$. The following conditions are
equivalent.

\begin{itemize}

\item [1.] The action $H\act\La_X H$ is 2-cocompact.

\medskip

\item [2.] $H$ is dynamically quasiconvex.

\end{itemize}
\end{thm}

\medskip

\noindent {\bf Remark.} Since every parabolic subgroup is always
dynamically quasiconvex we will  regard every group action on a
one-point set as 2-cocompact.

\subsection{Proof of the implication $2)\Rightarrow 1)$ }  The action $G\act X$ is
2-cocompact. So there exists a  compact fundamental set $K$ for
the action of $G$ on $\Th^2 X$. Denote by $\bu$ the entourage
$X^2\setminus K$. For every two distinct points $p$ and $q$ in $X$
there exists $g\in G$ such that $g(p,q)\in K$.   So $ (p,q)\not\in
\bu_1$ where $\bu_1=g^{-1}\bu.$ This means that the orbit $G\bu$
has separation property.

Let us first show that the index of $H$ in the stabilizer ${\rm
St}_G(\La H) =\{g\in G\ : g(\La H)=\La H\}$ of $\La H$ is finite.
Indeed for fixed two distinct points $\{p,q\}\subset\La H $ by the
exactness there exists   $\bu_1\in \U$ such that
$(p,q)\not\in\bu_1.$ So for every $h\in{\rm St}_G (\La H)$ we have
$ \{p,q\}\subset h(\La H)$ and $h(\La H)\not\in \s{Small}(\bu_1)$.
By the dynamical quasiconvexity applied to $\bu_1$ there are at
most finitely many such elements $h\in {\rm St}_G [\La H)$
distinct modulo $H.$ So ${\rm St}_X({\bf\La} H)=\cup_{j\in J}
k_jH\ (\vert J\vert <\infty).$

 Consider now the orbit $G(\La H)$. Applying
  the dynamical quasiconvexity  to $\bu$ we obtain a finite set
$\{g_i\in G\ :\ i\in I\} $ such that once $g(\bold\Lambda H)$ is
not $\bu$-small for some $g\in G$ then $g(\La H)=g_i(\La H)$ for
some $g_i.$

Consider   the following entourage of $\La H:$

$$\displaystyle \bv=\bigcap_{i\in I, j\in J} (k_j^{-1}g_i^{-1}\bu  \cap
\La^2H).$$

Let $(x,y)\in\Th^2({\bf\La} H)$. Since the orbit $G\bu$ has
separation property  there exists $g\in G$ such that
$g(x,y)\not\in\bu$ and hence $g(\La H)\not\in \s{Small}(\bu)$. We
have $g(\La H)=g_i({\bf\La} H)$ for some $i\in I.$ Hence
$g_i^{-1}g({\bf\La} H)=\La H$ and $g_i^{-1}g=k_jh$ for some $j\in
J$ and $h\in H.$ So $g(x,y)=g_ik_jh(x,y)\not\in \bu$. Consequently
$h(x,y)\not\in \bv$. We have proved that for every
$(x,y)\in\Th^2({\bf\La} H)$ there exists $h\in H$ such that
$(x,y)\not\in h^{-1}\bv$. This means that the set $\Th^2({\bf\La}
H)\setminus \bv$ is a compact fundamental set for the action
$H\act \Th^2(\La H)$. \bx

\subsection{Proof of the implication $1)\Rightarrow 2)$ }

We fix a group $G$ and a $32$-action of $G$ on a compactum $X$. By
Lemma \ref{cocomext} $G$ acts on the attractor sum $\ti X=X\cup\G$
where $\G$ is a connected, fine, hyperbolic graph and the action
$G\act\G^1$ is proper and cofinite. The canonical uniformity $\U$
of $\ti X$ is generated by an orbit $G\bu.$ By Lemma \ref{perdiv}
the restriction $\bu\vert_{\G^0}$ is a perspective divider.

\bigskip

Let $H<G$ be a subgroup. Denote by $\La H\subset\ti X$ its limit
set for the action on $\ti X.$  By Corollary \ref{hull} the set
$C=\HL$ is closed in $\ti X.$

\begin{lem}\label{degf} Every point $v$ in $ C^0=C\cap \G^0$ is either
  parabolic  for the action $H\act\La H$ or the number of edges
   incident to $v$ in the graph $C$ is
finite (i.e. the degree of $v$ in $C$ is finite).
\end{lem}

\proof Let $v\in C^0\setminus \La H$.  The set $\La H\subset X$ is
compact. So by the exactness of $\U$ there exists an entourage
$\bw\in \U$ such that $\{v\}\times \La H\cap\bw=\emptyset$. By the
visibility property there exists a finite set $F\subset \G^1$ such
that $\bu_F\subset \bw$. Hence   every eventual geodesic $\ga$
with endpoints $a=\ga_{-\infty}, b=\ga_{+\infty}\in\La H$ and
containing $v$ passes through $e_-$ and $e_+,$ where $e_-$ and
$e_+$ are edges of $F$ belonging to the geodesic rays joining
$\ga_{-\infty}$ with $v$ and $v$ with $\ga_{+\infty}$
respectively. Thus every
  arc $\ga$ in $C^1$ has a simple subarc $l$ between $e_-$ and $e_+$ which also
  contains $v$.
  Since $\ga$ has no intermediate stops $l$ is a geodesic.
   By the finess property of $\G$ the
number of such geodesic  subarcs   is finite.   Hence the number
of edges incident to $v$ is finite.

Suppose now that $v\in C^0\cap \La H.$ Then it is a parabolic
point for $G.$ Since $H$ acts 2-cocompactly on $\La H$ every point
of $\La H$ is either conical or bounded parabolic \cite[Main
theorem, b]{Ge1}. If $p$ is conical in $\La H$ then by
3-discontinuity of the action $G\act X$ it is also conical for
$G\act X$ which is impossible by \cite[Theorem 3A]{Tu2}. \bx

\begin{lem}\label{huldc} Let   $C=\HL.$ If $\vert C^1/H\vert <\infty$
then   $H$ is dynamically quasiconvex.

\end{lem}

\proof  We  extend the visibility property (see Definition
\ref{visprop}) from $\G^0$ to   $\ti X$. By Theorem \ref{exgeod}
for every $(x, y)\in\ti X^2$ there exists a geodesic $\ga\in\eG$
whose endpoints are $x$ and $y$. So for an edge $e\in \G^1$ and
$(x,y)\in \ti X^2$ put $(x,y)\in \bu_e$ if there exists such $\ga$
which does not contain the edge $e.$

Let $\bu\in \U$ and $\bv^3\subset \bu$.  Since the graph $\Gamma$
has the visibility property   there exists a finite set $E\subset
\G^1$ such that $\bu_E=\cap_{e\in E} \bu_e\subset\bv\vert_{\G^0}.$
Let $(x,y)\not\in\bu$ and let  $\ga$ be an eventual geodesic
joining $x$ with $y$. Choose $\{x', y'\}\subset \g$
  such that $(x,x')\in\bv$ and $(y,y')\in\bv$.  Since $(x,y)\not\in \bu$ we have
$(x',y')\not\in \bv\vert_{\G^0}.$ Hence $(x',y')\not\in\bu_E.$ So
  the piece of $\ga$ between  $x'$ and $y'$ contains an edge from $ E$.
 Hence $(x,y)\not\in\bu_e$. We have proved the inclusion $\bu_E\subset
\bu$ on  $\ti X^2$.

  If now $H$ is not dynamically
quasiconvex then by Definition \ref{dynquas}  the set $G_{\bold
u}=\{g\in G:g(\La H)\notin\s{Small}(\bold u)\}/H$ is infinite for
some $\bu\in\U.$ By the above argument there exists a finite
$E\subset \G^1$ such that $\bu_E\subset \bu$  on $\ti X.$ Since
$\vert \G^1/G\vert <\infty$  there exists an edge $e\in E$ for
which the set $\{g\in G : g(\La H)\not\in \s{Small}(\bu_e)\}/H$ is
infinite. Therefore the set $\{g\in G\ : \ e\in g (C^1)\}/H =\{
g\in G\ :\ g^{-1}(e)\in C^1\}/H$ is infinite too.  \bx

\medskip

 Suppose
that the action $H\act\La H$ is 2-cocompact. By Lemma \ref{huldc}
it is enough to prove that $\vert C^1/H\vert <\infty$ where
$C=\HL.$

Let $K$ be a compact fundamental set for the action $H\act
\Th^2(\La H)$. So $HK=\Theta^2(\La H)$. Let $\bu\in\U$ be an
entourage such that $\bu^3\cap K=\emptyset.$
 By
the visibility property there exists a finite set $F\subset\G^1$
such that $\bu_F\subset \bu.$ Thus $\bu_F^3\cap K=\emptyset.$ Up
to adding a finite number of edges to $F$ we can assume that $F$
is the edge set of a finite connected subgraph of $\G.$

 We call the edges of $C^1$ which belong to   $HF$
{\it red} edges. The other edges of $C^1$ are {\it white}.
Similarly we declare parabolic points of $H$   {\it red} and
other vertices of $C$ are {\it white.}

\begin{lem}\label{whiteray} Every infinite ray
$\rho:[0,\infty[\to C$ contains at least one red edge. Furthermore
every geodesic between two red vertices contains a red egde.

\end{lem}

\noindent {\it Proof of the lemma.}    By Lemma \ref{conv}   the
ray $\rho$ converges to a point $x=\rho(\infty)\in \La H.$ Since
the action $H\act\La H$ is 2-cocompact by \cite{Ge1} every point
of $\La H$ is either conical or bounded parabolic. By Proposition
\ref{geodray} $x$ is   conical for the action $H\act\La H.$

Therefore   there exists an infinite  set $S\subset H$ and two
distinct points $a$ and $b$ in $\La H$ such that $a$ and $b$ are
limit points for the sets $S(\rho(\infty))$ and   $S(\rho(0))$
respectively. Let $U_a$ and $U_b$ be disjoint $\bu$-small
neighborhoods of $a$ and $b$ for $\bu$ defined before the Lemma.
Thus $\exist s\in S\ :\ s(\rho(\infty))\in U_a,\ s(\rho(0))\in
U_b$. There exists $h\in H$ such that $h(a,b)\in K$. Hence $
h(a,b)\not\in \bu^3$. Since $h (s\rho(\infty),a)\in \bu$ and $h
(s\rho(0), b)\in \bu$ we obtain $
\partial(hs(\rho))\not\in\bu$. Thus $\partial(hs(\rho))\not\in\bu_F$
for the finite set $F$ defined therein.  It follows that $h
s(\rho)$ contains a red edge and so is $\rho.$

Let now  $\ga$ be a geodesic between two red points in $C$. Then
$\exist h\in H : h(\partial\ga)\in K$ so the pair $h(\partial\ga)$
is not $\bu_F$-small. Thus every geodesic $\ga$ connecting two red
points contains at least one red edge. The Lemma is proved. \bx

\medskip

 It remains to show that the set of white edges of $C^1$ is
$H$-finite.

Let us say that a segment of an eventual geodesic in $\C$ is {\it
 white} if all its edges and vertices are white. Denote by $\F$
 the  subgraph of $C^1$ obtained by adding to the set
 $F^1$ all adjacent  white segments. Since $F$ is connected
 $\F$ is also connected. By the
first statement of Lemma \ref{whiteray} every geodesic interval
containing only white edges
  has   finite length. Furthermore by Lemma \ref{degf} the
degree of every white vertex is finite. Thus by K\"onig Lemma the
connected subgraph $\F$ is finite.

We claim that $H\F^1=C^1.$ Indeed if $e=(a,b)\in\G^1$ is a white
edge then by the second statement of \ref{whiteray} one of its
vertices, say $a$, is white. Consider a maximal a white segment
$l_1$ of $C$ starting from $a$ and not containing $e$. It has a
finite length and ends either at a red vertex $c$ or at a red
edge. Our aim is to prove that the second case does happen for one
of such segments. Suppose it is not true for $l_1$. Then the other
vertex $b$ of the edge $e$ cannot be red. Indeed if $b$ is red,
then  $l_1\cup e$ has two red ends $c$ and $b$, and by Lemma
\ref{whiteray} $l$ must contains a red edge  which is impossible.
So $b$ is white. Then by Theorem \ref{exgeod} there exists another
maximal white segment $l_2$ starting from $b$. If it  ends up at a
red vertex $d$ then applying again \ref{whiteray} we obtain that
$l_1\cup l_2\cup e$ contains a red edge. So  there exists a  white
eventual geodesic segment $l$ starting from $e$ and terminating at
a red edge $e_1$. Thus there exists $h\in H\ :\ h(l\cup e)\subset
\F$. The Theorem is proved. \bx

\medskip

The proof of the above Theorem gives rise to another condition of
the dynamical quasiconvexity.

\begin{cor}
 \label{equiv} The following conditions are equivalent:

\begin{itemize}

\item[\sf a)] $H$ satisfies one of the conditions 1) or 2) of
Theorem \ref{infquas};

\medskip

\item[\sf b)] $\vert C^1/H\vert <\infty$ where $C=\HL.$
\end{itemize}
\end{cor}

\proof By Lemma \ref{huldc} it remains to prove that $\sf
a)\Rightarrow \sf b)$. By the statement $2)\Rightarrow 1)$ of
Theorem \ref{infquas} the dynamical quasiconvexity implies
2-cocompactness of the action $H\act\La_XH$. We have proved above
that the latter one implies that $\vert C^1/H\vert <\infty.$ \bx

\section{Pullback space for $32$-actions}
\label{secpulback}
 In  \cite[page 142]{Ge1} the following
 problem was formulated. Let a group $G$ admit
convergence actions on two compacta $T_i$  does there exist a
convergence action on a compactum $Z$ and two $G$-equivariant
continuous mappings $\pi_i : Z\to T_i\ (i=0,1)$  ?

\bigskip

\begin{center}
\begin{picture}(70,27)(-30,-20) \put(0,0){$Z$}
\put(-2,-2){\vector(-1,-1){20}}
\put(8,-2){\vector(1,-1){20}}\put(-25,-11){$\pi_0$}
\put(19,-11){$\pi_1$} \put(230,-20){(1)}\put(-30,-35){$T_0$}
\put(30,-35){$T_1$}\put(70,-35) {}
\end{picture}

\end{center}

\bigskip

\begin{dfn}\label{pullbackdef}
We call  the space $Z$  and  the action $G\act Z$ {\it pullback
space} and {\it pullback action} respectively.
\end{dfn}

Answering a question of
    M. Mitra \cite{M}  O.~Baker and T.~Riley   constructed in \cite{BR}
a hyperbolic group $G$ containing a free subgroup $H$ of rank $3$
such that the embedding does not induce an equivariant continuous
map (called ``Cannon-Thurston map'') $\partial  H\to
\partial G$ where $\partial $ is the boundary of a
hyperbolic group. Denote $T_0=\partial  H,$   and let
$T_1=\La_{\partial G} H$ be  the limit set for the action of
$H\act\partial G.$  The following proposition shows that
Baker-Riley's example is also a contre-example to the pullback
problem.

\begin{prop}\label{pullback}
The compacta $T_i\ (i=0,1)$ do not admit  a   pullback space on
which $H$ acts $3$-discontinuously.
\end{prop}

\proof Suppose by contradiction that   the diagram (1) exists.
Consider the spaces $\ti Z=Z\cup H,\ \ti T_0 = T_0\cup H,\ \ti
T_1=T_1\cup H$ equipped with the following topology (which we
illustrate only for $\ti T_0$ and  is defined similarly in the
other cases). A set $F$ is closed in $\ti T_0$ if

\medskip

\begin{itemize}

\item[1)] $F\cap T_0\in\s{Closed} (T_0);$

\medskip

\item[2)] $F\cap H\in \s{Closed}(H);$

\medskip

\item[3)] $\partial_1(F\cap H)\subset F$ where $\partial_1$
denotes the set of attractive limit points.

\end{itemize}

\medskip

The topology axioms are easily checked. Since $H$ is  discrete,
its points are isolated in $\ti T_0$ and the condition 2) is
automatically satisfied.

By \cite[Proposition 8.3.1]{Ge2} the actions $G\act \ti T_i$ and
$G\act \ti Z$ are 3-discontinuous. By the following lemma the maps
$\pi_i$ can be extended to the continuous maps $\ti\pi_0:\ti Z\to
\ti T_0$ and $\ti\pi_1:\ti Z\to \ti T_1$ where
$\ti\pi_i\vert_Z=\pi_i$ and $\ti\pi_i\vert_H={\rm id}\ (i=0,1).$

\medskip

\begin{lem}\label{extconv} Let $G$ be a group acting $3$-discontinuously
 on two compacta $X$ and $Y$.  Denote $\ti X$  and
 $\ti Y$ the spaces $X\cup G$ and $Y\cup G$ respectively  equipped with the above topologies.
 Suppose that the action on $Y$ is
 minimal and $\vert Y\vert > 2.$ If $f:X\to Y$ is a continuous $G$-equivariant
 map then the map   $\ti f:\ti
 X\to\ti Y$ such that $\ti f\vert_X=f$ and $\ti f\vert_G\equiv
 {\rm id}$ is continuous.
 \end{lem}

Assuming the lemma for the moment let us finish the argument.  By
hypothesis $H\act Z$ is 3-discontinuous. The map $\pi_0$ is
equivariant and continuous   and the action $H\act T_0$ is
minimal. So $\pi_0$ is surjective. Since $H$ is hyperbolic all
points of $T_0$ are conical \cite{Bo3}. By \cite[Proposition
7.5.2]{Ge2} the map $\pi_0$ is a homeomorphism. So we have the
equivariant continuous map $\pi=\pi_1\circ\pi_0^{-1} : T_0\to
T_1.$ By Lemma \ref{extconv} it extends equivariantly to the map
$\ti \pi: \ti T_0 \to\ti T_1$ where $\ti T_0=H\cup\partial H$ and
$\ti T_1=G\cup\partial_\infty G$. This is a Cannon-Thurston map. A
contradiction with the result of Baker-Riley. The Proposition is
proved modulo the following.

\medskip

 \noindent {\it Proof of  Lemma \ref{extconv}.} Let $F\subset \ti Y$
be a closed set.
 Denote $F_Y=F\cap Y$ and $F_G=F\cap G.$ We need to check that the set $\ti
 f^{-1}(F)=f^{-1}(F_Y)\cup F_G$ is closed. The conditions 1) and 2) are
 obvious for $\ti f^{-1}(F)\cap X$ and for $\ti
 f^{-1}(F)\cap G$ respectively.

  Let $z^\times=r\times X \cup X\times a$ be
 a limit cross for   $F_G$ on $X.$
 To check   condition 3) for the set $f^{-1}(F)$ we need to show that $b=f(a)\in F_Y.$
 Suppose not, and $b\not\in F_Y$
and let $B$ be a closed neighborhood of $b$ such that $B\cap
F_Y=\emptyset.$ Let $\bv\in\ent Y$ be an entourage such that
$B\bv\cap F_Y=\emptyset$ where $B\bv=\{y\in Y : (y,b_1)\in \bv,\
b_1\in B\}.$ Set $A=f^{-1}(B)\ni a.$ For a neighborhood $R$ of the
repelling point $r\in X$ the set $F_0=\{g\in F_G : g(X\setminus
R)\subset A\}$ is infinite.

Let $w^\times=p\times Y\cup Y\times q$ be a limit cross for $F_0$
on $Y,$ and $P\times Y\cup Y\times Q$ be its neighborhood. Since
$F_Y\subset Y$ is   closed  by condition 3)    we have $q\in F_Y$.
Suppose that $Q$ is $\bv$-small. By the hypothesis there exist
three distinct points
 $y_i\in Y\  (i=1,2,3)$.
Since the set $Y$ is minimal and $f$ is equivariant one has
$f^{-1}(y_i)=X_i\not=\emptyset$ and $X_i$ are mutually disjoint $
(i=1,2,3)$.

Let us now put some restrictions on $R$. Suppose   that $R\cap
X_i=\emptyset$ for at least two indices
$i\in\{k,j\}\subset\{1,2,3\}$ and
  for one of them, say $k$, we have $y_k\not\in P.$

   If
$g\in G$ is close to $w^\times$ we have $g(Y\setminus P)\subset Q$
and $g(y_k)\in Q.$ From the other hand $g(X_k)\subset A$  since
$X_k\cap R=\emptyset.$ Thus $g(y_k)\in Q\cap B$ and so
$(q,g(y_k))\in \bv.$ Hence $q\in B\bv$ and $q\not\in F_Y.$ A
contradiction. The lemma is proved.   \bx

\medskip

Since the answer to the   pullback problem for general convergence
actions is  negative, it  seems to be rather intriguing to study
the pullback problem in a  more restrictive case of 2-cocompact
actions. The rest of the section is devoted to a discussion of
this problem.

If   $G$ is a finitely generated group acting $3$-discontinuously
and $2$-cocompactly on compacta $X_1$ and $X_2$ then by the
Mapping theorem \cite[Proposition 3.4.6]{Ge2} there exist
equivariant maps $F_i :\partial G\to X_i\ (i=1,2)$ from the Floyd
boundary $\partial G$ of $G$. By \cite{Ka} the action on $\partial
G$ is 3-discontinuous. So $\partial G$ is a   pullback space for
any two $32$-actions of $G.$

If $G$ is not finitely generated  this argument  does not work as
the Mapping theorem requires the cofiniteness on edges of a graph
on which the group acts and which is not true for the Cayley
graphs
 in this case. An action of such a group
on a relative   graph    depends on the system of non-finitely
generated parabolic subgroups \cite[Proposition 3.43]{GePo2}.
Furthermore the action on the closure of the diagonal image of the
group in the product space may not be $3$-discontinuous. However
if  there is a pullback action for two 3-discontinuous actions of
$G$ and both of them are 2-cocompact then as shows the following
lemma a quotient of the pullback space also admits a $32$-action.

\medskip

\begin{lem}
\label{quotient} Suppose that $G$ acts 3-discontinuously and
2-cocompactly on two compacta $X_i\ (i=1,2)$. Let $X$ be a
pullback space for $X_i$ and $\pi_i:X\to X_i$ be  the
corresponding equivariant continuous maps. Then the action on the
quotient space
 $T=\pi(X)=\{(\pi_1(x,\pi_2(x)\ :\ x\in X\}$ is of type
 $(32)$.

Furthermore the action of $G$ on the attractor sum $\ti T=T\sqcup
G$ is also of type $(32)$.
\end{lem}

 \proof We will argue in terms of the attractor sums to obtain
 the more stronger last statement. By lemma \ref {extconv}
the maps $\pi_i$ extend to the continuous equivariant maps $\ti
\pi_i :\ti X\to\ti X_i$ where $\ti X = X\sqcup G$ and $\ti
X_i=X_i\sqcup G$, $\ti\pi_i\vert_G={\rm id},\ \ti\pi_i\vert
X_i=\pi_i$.

 By Lemma \ref{cocomext} the actions on
$X_i$ extends to 32-actions on the attractor sums $\ti X_i.$ Since
the action  $G\act \ti X_i$ is 2-cocompact there exists an
entourage  $\bu_i$ of $\ti X$ such that the uniformity $\U_i$ on
$\ti X_i$ is generated as a filter by the orbit $G\bu_i\ (i=1, 2)$
\cite[Proposition E, 7.1]{Ge1}.

Let $\ti\bu_i$ denotes the entourage ${\ti\pi_i}^{-1}(\bu_i)$ on
$\ti X (i=1, 2).$  Their $G$-orbits generate the lifted
uniformities $\ti \U_i.$ Then $\ti \bw=\ti\bu_0\cap \ti\bu_1$ is
an entourage of $\ti X$ whose orbit $G\ti\bw$ generates a
uniformity $\W$ on $\ti X$. Note that $\W$ is not a priori exact.
Indeed there could exists 2 points in $X$ such that
$\pi_i(x)=\pi_i(y)$ and there is no way to separate them using the
uniformities $\U_i\ (i=1, 2)$. So we consider the following
quotient spaces:

$$\ti T= \ti\pi(\ti X) =\{(\ti\pi_1( x), \ti\pi_2(x))\in \ti X_0\times\ti X_1\
:\  x\in \ti X\}, T=\pi(X)\ {\rm where}\ \pi=\ti\pi\vert_X. $$

Since $\ti\pi_i\ (i=0,1)$ are equivariant the map $\ti\pi$ is
equivariant too. Denoting by $\ti\pi_i:\ti T\to\ti X_{i-2}\
(i=3,4)$ the projections on the factors we obtain the following
commutative diagram.

\medskip

\begin{center}
\begin{picture}(70,27)(-30,-20) \put(0,0){$\ti X$}\put(2,-5){\vector(0,-1){50}}
\put(-5,-30){$\ti\pi$}\put(-1,-70){$\ti T$}
\put(-2,-2){\vector(-1,-1){30}}
\put(8,-2){\vector(1,-1){30}}\put(-30,-16){$\ti\pi_1$}
\put(25,-16){$\ti\pi_2$} \put(230,-20){(3)}\put(-45,-40){$\ti
X_1$} \put(40,-40){$\ti X_2$}
\put(-2,-62){\vector(-1,1){30}}\put(8,-62){\vector(1,1){30}}
\put(-18,-45){$\ti\pi_3$}\put(12,-45){$\ti\pi_4$}
\end{picture}

\end{center}

\vspace*{2cm}

  Since the map  $\pi$ is continuous
and surjective the action $G\act \ti T$ is 3-discontinuous  too
\cite[Proposition 3.1]{GePo1}. It remains to prove that it is
2-cocompact in the quotient topology. Let $\ti\bv=\pi(\ti\bw)$ and
consider the uniformity $\ti\V$ on $\ti T$ generated by the orbit
$G\ti\bv$. To show that $G\act\ti T$ is 2-cocompact it is enough
to prove that $\ti\V$ is exact \cite[Proposion E, 7.1]{Ge1}.

So let $x, y$ be two distinct points of $\ti T$. Then  either $\ti
\pi_3(x)\not= \pi_3(y)$ or $\ti\pi_4(x)\not=\ti\pi_4(y)$. For
example in the first case by the exactness of ${\ti\U}_0$ there
exists $g\in G$ such that $g(\ti \pi_3(x),\ti\pi_3(y))\not\in
\ti\bu_0$. By definition of $\ti \bw$ we obtain that
$(\ti\pi(x),\pi(y))\not\in \ti\bv.$ So $\ti\V$ is exact. \bx

\medskip

The aim of the following Proposition is to show that two
$32$-actions may not have a pullback.  We note that it is one of
the rare cases when a fact known for   finitely generated
relatively hyperbolic groups   is not in general true for
non-finitely generated groups.

\begin{prop}\label{freeinf}
The free group $F_{\infty}$ of countable rank admits two
$32$-actions not having a   pullback space.
\end{prop}

\proof Let $G=<x_1,...,x_n, y_1,...,y_m,...>\ (n\geq 2)$ be a
group freely generated by the union of a  finite set
$X=\{x_1,...,x_n\}$ and an infinite set $Y=\{y_1,...\}$. Let
$A=<X>$ be a subgroup generated by $X$, and let  $H$ be a subgroup
of $A$ freely generated by an infinite set $W=\{w_i\ :\ i\in\N\}.$

Set $Z=\{z_m=y_m w_m\ :\ m\in\N\}$, $P= <Y>$ and $Q=<Z>.$ The set
$X\cup Z$ can be obtained by Nielsen transformations from $X\cup
Y$ \cite{LS}. So  $X\cup Z$ is also a free basis for $G,$ and the
map $\varphi\ :\ x_i\to x_i, y_k\to z_k \ (i=1,...,n; k\in\N)$
  extends to an automorphism of $G$.  We have  two splittings of $G:$

 $$G=A*P,\ \ {\rm and}\ \ G=A*Q.\hfill\eqno(1)$$

Each splittings in (1)   gives rise to an action of $G$ on a
simplicial tree   whose vertex groups are conjugates of either $A$
or  $P$ (respectively to $Q$). We now replace the vertices
stabilized by $A$ and its conjugates by their Cayley trees. Denote
the obtained simplicial $G$-trees by $\T_i\ (i=1,2)$. Their edge
stabilizers are trivial and   vertex stabilizers are non-trivial
if only if they are conjugate to $P$ (respectively to $Q$). The
vertices of  $\T_1$ (respectively $\T_2$) are   the elements of
$G$ and the parabolic vertices corresponding to conjugates of $P$
(respectively $Q$). The graph $\T_i$ is a connected fine
hyperbolic graph such that the action of $G$ on edges are proper
and cofinite. Hence the actions satisfy Bowditch's criterion of
relative hyperbolicity \cite{Bo1}. By \cite{Ge2} both actions on
the trees extend to $32$-actions on compacta $R_i$ which are the
limit sets for the actions $G\act\T_i\ (i=1,2)$.

We  claim that $P\cap g^{-1}Qg=\{1\}$ for all $g\in G.$ Indeed
consider the endomorphism $f$ such that $f(x_i)=x_i,\ f(y_j)= 1\
(i=1,...,n, j=1,...)$. The map $f$   is injective on $Q$ as well
as on every conjugate   $g^{-1}Qg$. From the other hand $Y\subset
{\rm Ker} f$. So $P< {\rm Ker} f.$ We have proved that
$$\foral  g\in G : P\cap g^{-1}Qg=\{1\}.\hfill\eqno(2)$$

Arguing now by contradiction assume that there exists a pullback
space $R$ and equivariant projections $\pi_i : R\to R_i\ (i=1,2)$.
By Lemma \ref{quotient} the action on the quotient space:

$$ T= \pi( R) =\{(\pi_1( r), \pi_2(r))\
\vert\  r\in R\}
$$

\noindent is 3-discontinuous and 2-cocompact. Note that the action
$G\act T$ is minimal because $G\act R$ is minimal.

By \cite[Main theorem, b]{Ge1} all points of $T$ are either
conical or bounded parabolic. If  $p\in T$ is parabolic  then
$\pi_{i+2}(p)$  are parabolic points in both $ R_i$ for the map
$\pi_{i+2}=\ti\pi_{i+2}\vert_T$ (see the diagram in Lemma
\ref{quotient}). Indeed the preimage of a conical point by an
equivariant map is conical \cite[Proposition 7.5.2]{Ge2}. So $p$
must be fixed by the intersection of some parabolic subgroup
$g_1Pg_1^{-1}$ of the first action and a parabolic subgroup
$g_2Qg_2^{-1}$ of the second
 ($g_i\in G$). However by (2) this intersection is  trivial. Thus there
  are no parabolic points for the $32$-action $G\act
 T$. By \cite[Theorem 8.1]{Bo3} (see also \cite[Corollary
3.40]{GePo2}) the group $G$ is hyperbolic and so finitely
generated. This is a contradiction.\bx

\bigskip

\noindent {\bf Definition.}   The set of the stabilizers of the
parabolic points   for a $32$-action on a compactum $X$ is called
{\it peripheral structure} on $G.$

\medskip

 The
following theorem  provides a sufficient condition for the
existence of pullback space for two $32$-actions of a group.

\medskip
\begin{theor}\label{suffcond} Let $G$ be a group which admits $32$-actions on compacta
$X$ and $Y$. Let $\P$   be the peripheral structure corresponding
to the action on $X$. Suppose that every $P\in \P$ acts
2-cocompactly on $\La_YP$. Then there exists a compactum $Z$
equipped with a $32$-pullback action of $G$ with respect to its
actions on $X$ and $Y$.
\end{theor}

\medskip

\medskip

\noindent {\bf Remark.} Using Theorem \ref{infquas} one can
reformulate the hypotheses above by requiring that the action of
each subgroup $P\in\P$    is dynamically quasiconvex on their
limit sets in $Y$.

\medskip

 \noindent {\it Proof of Theorem \ref{suffcond}.} Denote
$$\mathcal R=\{P \cap Q:P\in\mathcal
P,Q\in\mathcal Q,|P \cap Q|=\infty\}.\hfill\eqno(4)$$

We will indicate a compactum $Z$ acted upon by $G$
$3$-discontinuously and $2$-cocompactly whose peripheral structure
is $\Rc.$ Denote by $\Par(Y,P)$ the set of parabolic points for
the action of $P\in\P$ on $Y.$ We will need the following lemma.
\begin{lem}
\label{parabcor} Let $G, X,Y, \P, \Q$ be as in Theorem B. The
following properties of a subgroup $H{\subset}G$ are equivalent:
\begin{itemize}

\item [a:] $H\in\mathcal R$;

\item [b:]  there exist $P\in\mathcal P$ and $q\in\Par(Y,P)$ such
that $H=\mathsf{St}_Pq$.
\end{itemize}
\end{lem}
 \noindent {\it Proof of the Lemma.}
$\mathsf{b{\Rightarrow}a})$. If the action of $P\in\P$ on $Y$
admits a parabolic point $q$ then its stabiliser $H$ is an
infinite subgroup of $P$. By \cite[Theorem 3.A]{Tu2} the point $q$
is parabolic
 for the action $G\act Y$. We note that the assumption of \cite{Tu2}
 that the space is metrisable can be omitted by a small
 modification of the argument. Let $Q={\rm St}_Gq\in\Q$ be the stabilizer of
$q$. We obtain $H=P\cap Q\in\Rc.$

{\it $\mathsf{a{\Rightarrow}b})$.} Let $H=P \cap Q$ for
$(P,Q)\in\mathcal{P{\times}Q}$. We may assume that
$P=\mathsf{St}_Gp$ and $Q=\mathsf{St}_Gq$ for $p\in\Par(X,G)$,
$q\in\Par(Y,G)$. Since $H$ is an infinite subset of $Q$ we have
$\La_YH=\{q\}$. Since $H{\subset}P$ we have $q\in\La_YP$. If $q$
is conical for $P{\curvearrowright}\La_YP$ then it is also conical
for $G{\curvearrowright}Y$   contradicting by \cite[Theorem
3.A]{Tu2} the fact that $q\in\Par(Y,G)$. The lemma is proved.

\medskip

 The peripheral structure  $\P$ consists of finitely
many $G$-conjugacy classes \cite[Main Theorem, $\mathfrak
a$]{Ge1}. Since for every $P\in\P$ the action $P\act\La_YP$ is
2-cocompact   there are   finitely many $P$-conjugacy classes of
maximal parabolic subgroups in $P.$  So it follows from Lemma
\ref{parabcor} that $\Rc$ consists of finitely many $G$-conjugacy
classes. Since the subgroups in $\mathcal R$ are infinite, each of
them is contained in exactly one $P\in\mathcal P$ and in exactly
one $Q\in\mathcal Q$. So the inclusions induce well-defined maps
$\mathcal P\overset\pi\longleftarrow\mathcal
R\overset\sigma\longrightarrow\mathcal Q$ equivariant by
conjugation.

We now extend the maps $\pi,\sigma$  identically over the sets
$\widetilde \P=G{\sqcup}\mathcal P$, $\widetilde
\Q=G{\sqcup}\mathcal Q$, $\widetilde R=G{\sqcup}\mathcal R$.
Denoting the extensions by the same symbols  we have
$G$-equivariant maps $\widetilde{\mathcal
P}\overset\pi\longleftarrow\widetilde{\mathcal
R}\overset\sigma\longrightarrow\widetilde{\mathcal Q}$. By Lemma
\ref{cocomext} the set $\ti P$ is the vertex set of a connected
  fine graph $\De$ such that the action on edges
$G\act\De^1$ is cofinite and proper.

We will  construct a   connected graph $\G$ whose vertex set is
$\G^0=\ti \Rc$   and the action on edges is cofinite and proper.
  The set of edges $\G^1$  will be  obtained by replacing the
parabolic vertices $\P$ of $\De$ by connected graphs  coming from
the  action of the groups $P\in\P$ on $Y.$ We do it in the
following four steps.

\medskip

\noindent {\bf Step 1.}  {\it Definition of $\Gamma^1_1$.}

 Choose a set $\mathcal R_\#\subset \mathcal R$ that intersects each
conjugacy class by a single element. For every $R\in\mathcal R_\#$
we join the vertex $R$ with each element of $R\subset G$ and
denote by $E_R$ this set of edges.  Then put

$$\displaystyle \G^1_1=\bigcup  \left\{ gE_R\ : g\in G,\ R\in R_\#\right\}.$$

The set $\G^1_1$ corresponds to the well-known coned-off
construction over every coset $gR$ where $g\in G, R\in \Rc_\#$
\cite{Fa}, \cite{Bo1}.

\medskip

\noindent {\bf Step 2.} \noindent {\it Definitions of $\G_2^1$ and
$\G_3^1$.}

Choose a set $\mathcal P_\#\subset \mathcal P$ that intersects
each conjugacy class of $P$ by a single element. For each
$P\in\mathcal P_\#$ we add to $\Gamma^1$ a connected $G$-finite
set of pairs according to one of the following ways.

\medskip

\noindent {\bf Case 2.1.} ({\it hyperbolic case}) \noindent
$P\in\mathcal P_\#{\setminus}\mathsf{Im}\pi$ (or equivalently
$\pi^{-1}(P)\cap\Rc=\emptyset$).

Then $P$ acts on $Y$ either as an elementary loxodromic $2$-ended
subgroup, or the action $P\act \La_YP$ is a non-elementary
32-action without parabolics \cite[Theorem 3A]{Tu1}. In both cases
every point of $\La_YP$ is conical \cite[Main Theorem, $\mathfrak
b$]{Ge1} and $G$ is a hyperbolic group \cite [Theorem 8.1]{Bo3}
 (for another proof of this fact
see \cite[Appendix]{GePo1}).  There exists a $P$-finite set
$\Gamma^1_P$ of pairs of elements of $P$ such that the graph
$(P,\Gamma^1_P)$ is connected.

\medskip

Put $\displaystyle \G_2^1 =\bigcup\left\{ g\G^1_P\ :\  g\in G,\
P\in\mathcal P_\#{\setminus}\mathsf{Im}\pi\right\}.$

\medskip

\noindent {\bf  Case 2.2} ({\it non-hyperbolic case})
$P\in\mathcal P_\# \cap \mathsf{Im}\pi$.

\noindent There is a canonical bijection
$\tau_P:\Par(Y,P)\to\pi^{-1}P$. Let

\medskip

\centerline{$M_P=\{g \in G:g$ is joined by a $\Gamma_1^1$-edge
with some $R\in\pi^{-1}P\}$.}

\medskip

The set $M_P$ is $P$-invariant and $P$-finite by the construction.
Let $\G_P^1$ be the $P$-finite set of pairs of the elements of the
$P$-invariant set $M_P{\cup}\pi^{-1}P$ such that $(P, \G^1_P)$ is
connected. The latter one exists as the graph corresponding to the
$32$-action $P\act\La_YP$ is connected by Lemma \ref{cocomext}.
Put

$$\G_3^1=\bigcup\left\{g\G_P^1 \ :\ g\in G,\ P\in\mathcal
P_\# \cap \mathsf{Im}\pi\right\}.$$

\medskip

\noindent {\bf Step 3.} {\it Definition of $\G_4^1$.} Consider in
the graph $\De$ the set of all its ``horospherical'' edges
$\De_0^1=\{(P,g)\ :\ P\in\P, g\in M_P\}$. Let $\G_4^1$ be the set
of all non-horospherical edges of $\De$:

$$\G_4^1=\De^1\setminus \De^1_0.$$

\medskip

\noindent {\bf Step 4.}  {\it Definition of $\G^1.$} Let
$$\G^1=\G_1^1\cup\G_2^1\cup\G_3^1\cup\G_4^1.\hfill\eqno(5)$$

The set $\G^1$ is obviously $G$-finite.

\begin{lem}\label{connect} The graph $\Gamma=(\G^0, \G^1)$ is connected.
\end{lem}

\noindent {\it Proof of the lemma.}  Since every vertex of
$\Gamma$ is either an element of $G$ or is joined with an element
of $G$ it suffices to verify that every two elements of $G$ can be
joined by a path in $\Gamma$. We initially join them by a path
$\gamma$ in the connected graph $\Delta$. We transform this path
as follows. If all vertices of $\gamma$ belong to $G$ then
$\gamma$ is also a path in $\Gamma$. It $\gamma$ passes through a
point $P\in\mathcal P$ then it has a subpath of the form
$g_0{-}P{-}g_1$ where $g_0,g_1 \in M_P$. The graph of $P$
corresponding to the action on $Y$ is connected. So we can replace
this subpath by a subpath with the same endpoints all whose
vertices are contained in $M_P$ (in the ``hyperbolic'' case) or in
$M_P{\cup}\pi^{-1}P$ (in the ``non-hyperbolic'' case). In both
cases the edges of this new subpath belong to $\Gamma^1$. The
lemma is proved.

\medskip

\noindent {\it End of the proof of Theorem \ref{suffcond}.} By
Lemma \ref{perdiv}  the sets $\ti \P$ and $\ti Q$ admit
perspective dividers $\mathbf u\subset \widetilde{\mathcal P}^2$,
$\mathbf v\subset \widetilde{\mathcal Q}^2$. Since the projection
$\pi$ and $\si$ commute with the group action  the lifts
$\pi^{-1}\mathbf u$ and $\sigma^{-1}\mathbf v$ are perspective
dividers on $\G^0=\widetilde{\mathcal R}$.

It is   a direct verification that $\bw=\pi^{-1}\mathbf
u\cap\sigma^{-1}\mathbf v$ is a perspective divider on $\G^0$.
Indeed if $g(a,b)\not\in\bw^0$ then $g(\pi(a), \pi(b))\not\in
(\bu^0=\bu\vert_{\ti\P})$ or $g(\si(a), \si(b))\not\in
(\bv^0=\bv\vert_{\ti\Q})$. So there exist at most finitely many
such elements $g\in G$ as $\bu^0$ and $\bv^0$ are both
perspective. Similarly $\bw^0$ is a divider on $\G^0$ as if $(\cap
F_1\{\bu^0\})^2\subset \bu^0$ and $(\cap F_2\{\bv^0\})^2\subset
\bv^0$ for some finite $F_i\subset G\ (i=1,2)$ then $(\cap
F\{\bw^0\})^2\subset \bw^0$ where $F=F_1\cap F_2.$

It follows that the projections $\pi:(\G^0, \bw)\to (\ti \P, \bu)$
and $\si:(\G^0, \bw)\to (\ti \Q, \bv)$ are uniformly continuous
with respect to the uniformities generated by the divider orbits.

  By Lemma \ref{compldense} the action of $G$ on the
Cauchy-Samuel completion $\ti Z$ of $(\G^0, \bw)$ is a
$32$-action. By \cite[II.23, Proposition 13]{Bourb} the completion
$\ti Z$  coincides with the closure ${\rm Cl}_{\ti X\times \ti
Y}(\G^0)$
  of $\G^0$ embedded diagonally in $\ti X\times \ti Y$ where $\ti X=X\sqcup G$
  and $\ti Y=Y\sqcup G.$
So the projections $\pi$ and $\si$ extend continuously to the
equivariant maps $\ti\pi:\ti Z\to\ti X$  and $\ti\si:\ti Z\to \ti
Y$ whose restrictions to $G$ is the identity. We have proved that
$Z=\La_{\ti Z}G$ is a pullback space.  The Theorem is proved. \bx

\bigskip

\noindent To prove the  statement converse to Theorem
\ref{suffcond} we need the following direct generalization  of the
argument of \cite[Lemma 2.3, (4)]{MOY}  avoiding the metrisability
assumption.

\begin{lem}\label{preimlim} Let a group $G$   admits two
non-trivial 3-discontinuous  actions on compacta $X$ and $Y$, and
let   $f:X\to Y$ be an equivariant continuous map. Let $H$ be a
subgroup of $G$ such that $\La_YH\subsetneqq Y$. Suppose that $H$
acts cocompactly on $Y\setminus\La_YH.$ Suppose that  for every
infinite set $B\subset G\setminus H$ there exist an infinite
subset $B_0\subset B$  and at least two distinct points $r_i\in
f^{-1}(\La_YH)$ such that $\foral  g\in B_0 : g(r_i)\not\in
f^{-1}(\La_YH)\ (i=1,2).$ Then $f^{-1}(\La_YH)=\La_XH.$
\end{lem}

\medskip

\begin{cor}\label{preimp}
If $p$ is a bounded parabolic point for the action of $G$ on $Y$
then $f^{-1}(p)$ is the limit set $\La_X({\rm St}_Yp)$ of ${\rm
St}_Gp$ for the action on $X.$
\end{cor}

\medskip

\noindent {\it Proof of the Corollary.} By the equivariance and
continuity of $f$ we have $\La_XH\subset f^{-1}(\La_YH)$ where
$H={\rm St}_Gp$. So if $f^{-1}(p)$ is a single point then the
statement is trivially true. If $f^{-1}(p)$ contains at least two
distinct points $r_i\ (i=1,2)$ then we have $\foral  g\in
G\setminus H\ :\ g(r_i)\not\in f^{-1}(\La_YH)$  as $g(p)\not=p.$
and we apply  Lemma \ref{preimlim}.\bx

\medskip

{\it Proof of the Lemma.} The statement is trivial if $H$ is
finite, so we assume that $H$ is infinite. Suppose first that the
set $f^{-1}(\La_YH)$ is   finite. Since $f(\La_XH)\subset \La_YH$
then $f^{-1}(\La_YH)$ is pointwise fixed under a finite index
subgroup of $H$. So   $f^{-1}(\La_YH)=\La_XH$ in this case.

Suppose that $f^{-1}(\La_YH)$ is infinite.  Suppose by
contradiction that there exists a point $s\in
f^{-1}(\La_YH)\setminus \La_XH$. Then there exist an infinite set
$B\subset G\setminus  H$ converging to the cross whose attractive
limit point is $s$. By our assumption  there exists an infinite
subset $B_0\subset B$  and distinct points $r_i\in f^{-1}(\La_YH)$
such that $\foral g\in B_0\ : g(r_i)\not\in f^{-1}(\La_YH)\
(i=1,2).$ Then   one of them  $z\in\{r_1, r_2\}$ is not repulsive
for the limit cross of $B_0.$ So   for every open neighborhood
$U_s$ of $s$ there exists an  infinite subset
 $B'_0\subset B_0$ such that $\foral g\in
B'_0\ :\ g(z)\in U_s\setminus f^{-1}(\La_YH)$.

Let $K$ be a compact fundamental set for the action $H\act
(Y\setminus\La_YH).$ Since $X$ is compact and $f$ is equivariant
the set $f^{-1}(K)=K_1$ is a compact fundamental set for the
action of $H$ on $X\setminus f^{-1}(\La_YH)$. Therefore  for every
$g\in B$ there exists $h\in H$ such that $hg(z)\in K_1.$ The set
$$A_s=\{h\in H : h(K_1)\cap U_s\not=\emptyset\}$$ is infinite for
every open neighborhood $U_s.$ Indeed if it is not true for some
$U_s$   then by the argument above the orbit $A_s(K_1)$ intersects
every  neighborhood $U^*_s$ of $s$ such that $U^*_s\subset U_s.$
Then by compactness of $K_1$ we would have $h^{-1}(s)\in K_1$ for
some $h\in H,$ implying that $f(s)\in h(K).$ This is impossible as
$\La_Y(H)\cap h(K)=\emptyset$ for any $h\in H.$ Therefore there
exists  infinitely many $h\in H$ such that $h(K_1)\cap
U_s\not=\emptyset$ for every neighborhood $U_s$ of $s$. Thus $s\in
\La_XH.$ A contradiction. \bx

\bigskip

\noindent The main result  of the paper is the following.

\begin{thm}\label{crit}

 Two $32$-actions of $G$ on
compacta $X$ and $Y$ with peripheral structures $\P$ and $\Q$
admit a pullback space $Z$ if and only if one of the following two
conditions is satisfied (and so both of them):

\begin{itemize}
\item [1.]   $\foral P\in \P$ acts 2-cocompactly on $\La_YP$
\item[2.] $\foral Q\in\Q$ acts 2-cocompactly on $\La_XQ.$
\end{itemize}

\medskip
\end{thm}

 \proof After Theorem \ref{suffcond} we need only to show that if the pullback space $Z$
 exists then every $Q\in\Q$ acts
 2-cocompactly on $\La_XQ.$ Suppose that $G$ admit a pullback action
  $G\act Z$ for two
$32$-actions on $X$ and $Y$ and let $\displaystyle X\overset{f_1}
{\leftarrow} Z\overset{f_2} {\rightarrow}   Y$ be the equivariant
continuous maps.

By Lemma \ref{quotient}   we may assume that the action  $G\act Z$
is 2-cocompact. Let $q\in Y$ be a  parabolic  point and $Q={\rm
St}_Gq\in \Q$.

\medskip

  By Corollary \ref{preimp}
$f_2^{-1}(q)$ is the limit set $\La_ZQ$. Since
$(Y\setminus\{q\})/Q$ is compact, $Z$ is compact and $f_2$ is
equivariant and continuous, $Q$ acts cocompactly on $Z\setminus
f_2^{-1}(q)=Z\setminus\La_ZQ$.

The set $f_1(\La_ZQ)$ is a closed $Q$-invariant subset of $X$. So
$\La_XQ\subset f_1(\La_ZQ).$ Since $f_1$ is continuous and
equivariant   we have $\La_XQ= f_1(\La_ZQ)$ and the action $Q\act
X\setminus \La_XQ$ is cocompact.

By Lemma \ref{cocomext} there exists a connected, fine, hyperbolic
graph $\G_1$ corresponding to the 32-action $G\act X$ such that
the action $G\act \G^1_1$ is proper and cofinite.

In the following lemma we will  use the notion of a dynamical
bounded subgroup introduced in \cite[section 9.1]{GePo3}. Recall
the topological version of this definition: a subgroup $Q$ of $G$
is said to be {\it dynamically bounded} for the action on $X$ if
there exist finitely many proper closed subsets $F_i$ of $X$ such
that $\foral g\in G\ \exist i\ :\ g(\La_XQ)\subset F_i$.

Considering the action of $Q$ on $\ti X= X\cup\G_1$ we have.

\begin{lem}\label{dynbound}
If a subgroup $Q$ of $G$ acts cocompactly on $X\setminus\La_XQ$
then it acts cocompactly on $\ti X\setminus\La_XQ$.
\end{lem}

\noindent {\it Proof of the lemma.} The parabolic subgroup $Q$ is
obviously dynamically bounded for the action on $Y$. Indeed since
$Y$ is compact there are finitely many closed proper subsets $R_i$
of $Y$ such that $\foral g\in G\   \exist i\in\{1,...,m\}\ :\
g(q)\in R_i.$   Since $F_i=f_2 f_1^{-1}(R_i)$ is closed in $X$ and
$f_i$ is surjective and equivariant, we have $\displaystyle
X=\bigcup_{i\in\{1,...,m\}} F_i$ and $g(\La_XQ)=g(f_2
f_1^{-1}(\La_YQ))\subset F_i$ for some $i\in\{1,...,m\}.$ So $Q$
is dynamically bounded for the action on $X.$

The proof of  \cite[Proposition 9.1.3]{GePo3} implies that if $Q$
is dynamically bounded on $X$ and   acts cocompactly on $
X\setminus \La_XQ$ then it acts cocompactly on $\ti X\setminus
\La_XQ$.\bx

\medskip

\noindent {\bf Remark.} The metrisability assumption stated in
\cite[Proposition 9.1.3]{GePo3}
 was only used   to satisfy
another (metric) definition of the dynamical boundness which we do
not use here.

\medskip

{\it End of the proof of Theorem \ref{crit}.} By Lemma
\ref{dynbound} the action $Q\act (\ti X\setminus\La_XQ) $ is
cocompact. Let $K\subset(\ti X\setminus\La_XQ) $ be a compact
fundamental set for this action. By Corollary  \ref{equiv} it is
enough to prove that $\vert C^1/Q\vert <\infty$ where $C^1$ is the
set of edges of $C={\rm Hull}_X(\La Q).$ Let $e=(a,b)\in C^1.$
Then one of its vertices, say $a$, is not in $\La_XQ$.  By
definition of $C$ there exists an infinite eventual geodesic $\g$
such that $e\subset \g(\Z)$ and $\g(\{-\infty,
+\infty\}\subset\La_XQ.$ So there exists $g\in Q$ such that
$g(a)\in K\cap C$ and $ge\subset g\ga(\Z).$

  We have
$\La_XQ\cap K=\emptyset$.  By the exactness of the uniformity $\U$
of the topology $\ti X$ there exists an entourage $\bu\in\U$ such
that $\bu\cap (K\times\La_XQ)=\emptyset.$ By the visibility
property there exists a finite set $F\subset \G^1$ such that
$\bu_F\subset \bu.$ So every geodesic from $K$ to $\La_XQ$
contains an edge from $F.$ Hence $g\ga(\Z)$ contains a finite
simple geodesic subarc $l$ such that $g(a)\in l^0$ and  $\partial
l\subset F^0.$ Since that the graph $\G$ is fine there are
finitely many geodesic simple arcs joining the vertices of $F^0$.
So the set $E$ of the edges of these arcs is finite. We have
proved that $Q(E)=C^1$ and so $\vert C^1/Q\vert<\infty$. Theorem
\ref{crit} is proved. \bx

\section{Corollaries}
\label{equivsec}

\noindent The goal of this Section is the following list of
corollaries.

\medskip

\begin{cor}\label{equivac}  Let a group $G$ acts on compacta
$X$ $3$-discontinuously and $2$-cocompactly. Let $\P$ and $\Q$ be
the peripheral structures   for the actions on $X$ and $Y$
respectively. Then the following statements are true.

\begin{itemize}

\item [a)] Suppose that one of the conditions 1) or 2) of Theorem
\ref{suffcond} is satisfied then    $G$ is relatively hyperbolic
with respect to the system $\Rc=\{P\cap Q\ :\ P\in\P, Q\in \Q,
\vert P\cap Q\vert=\infty\}.$

\medskip

\item [b)] Every $P\in \P$ acts $2$-cocompactly on $\La_YP$ if and
only if every $Q\in\Q$ acts 2-cocompactly on $\La_XP$.

\medskip

\item [c)] Assume that $\foral P\in \P\ \exist Q\in\Q : P<Q.$ Then
the induced action of every $Q\in\Q$ on $\La_XQ$ is 2-cocompact.

\medskip

\item [d)] Suppose that for every $P\in\P$ there exists $Q\in\Q$
such that $P < Q.$ Then there exists an equivariant continuous map
$f: X\to Y $. Furthermore the induced action of every $Q\in\Q$ on
$\La_XQ$ is 2-cocompact.

\medskip

 \item [e)] If $\P=\Q$ then the actions are equivariantly
 homeomorphic.
\medskip

\item [f)] Let $G$ admits a 32-action on a compactum $X$ and $H<
G$ is a parabolic subgroup for this action. Then we have:

\begin{itemize}

\item[f1)] If $G$ is finitely generated then for any other
32-action $G\act Y$ the subgroup $H$ is dynamically quasiconvex.

\item[f2)] If $G$ is not finitely generated the statement f1) is
not true in general.
\end{itemize}

 \end{itemize}

\end{cor}

\medskip

\proof  a) directly follows from Theorem \ref{suffcond}.

\medskip

b) By Theorem \ref{suffcond} the pullback space $Z$ exists. Then
by Theorem \ref{crit} every $Q\in\Q$ acts 2-cocompactly  on
$\La_XQ$.

\medskip

c) By the assumptions the elements of $\P$ act parabolically on
$Y.$ So they  all act 2-cocompactly on their limit sets on  $Y$.
Then by b) the elements of $\Q$ act 2-cocompactly on their limit
sets on $X.$

d) By the assumptions the elements of $\P$ act parabolically on
$Y$. So they act 2-cocompactly on their limit sets on $Y$. By
Theorem \ref{suffcond} there exists a $32$-action $G\act Z$ which
is a pullback action for the actions on $X$ and $Y$. We have two
equivariant continuous maps $\pi : Z\to X$ and $\si :Z\to Y.$

We claim that  $\pi$ is injective. Indeed every point $x$ of $X$
is either conical or bounded parabolic \cite[Main Theorem,
b]{Ge1}. If     $x\in X$ is conical then    $\pi^{-1}(x)$ contains
is a single point \cite[Proposition 7.5.2]{Ge2}.

If $p\in X$ is a bounded parabolic  and $P={\rm St}_G p$ then by
Corollary \ref{preimp}  $\pi^{-1}(p)=\La_Z(P)$. By   Theorem
\ref{suffcond}  the peripheral structure  for the action $G\act Z$
is   $\Rc=\P\cap \Q=\P.$ By  the assumption we have $P\in \Rc$ is
  parabolic   for the action on $Z$, so $\pi^{-1}(p)$ is
a single point. Thus   $\pi$ is injective and so is a
homeomorphism. Hence map $  f=\si\circ\pi^{-1}:X\to Y$  is
equivariant and continuous.

\medskip

e) follows from d).

\medskip

f1) Indeed if $G$ is finitely generated the Floyd boundary
$\partial_fG$ is universal pullback space for any two $32$-actions
of $G$ on $X$ and $Y$ \cite[Map theorem]{Ge1}. So if $H$ is a
maximal parabolic subgroup for the action on $X$ by the necessary
condition  of Theorem B it acts 2-cocompactly on $\La_YH.$ Then by
Theorem A it is dynamically quasiconvex.

\medskip

f2) Proposition \ref{freeinf} provides a contre-example. Indeed
there is no pullback action for two 32-actions of the free group
$F_\infty$ on two spaces. By Theorem B  there exists a parabolic
subgroup   of one of the actions  which does not act 2-cocompactly
on its limit set for the other one. Again by Theorem A the sugroup
is not dynamically quasiconvex for the second action.

The Corollary is proved. \bx

\medskip

\begin{rems}\label{conapar}  {\rm  The statements d) and e) give rise to
more restrictive similarity properties of $32$-actions given by
equivariant maps.

The statement d) was already known in several partial cases. If
first, $G$ is finitely generated then it follows from the
universality of the Floyd boundary. Indeed by \cite{Ge2} there
exist continuous equivariant (Floyd) maps $F_1:\partial G\to X$
and $F_2:\partial G\to Y$ where $\partial G$ is the Floyd boundary
of the Cayley graph of $G$ (with respect to some admissible scalar
function). By \cite[Theorem A] {GePo1} for a parabolic point $p\in
\La_X G$ the set $F_1^{-1}(p)$
  is the limit set $\La_{\partial G} P$ of the stabilizer $P={\rm St}_{G}p$ for
   the
action $G\act\partial G$. Since  $F_2$ is equivariant the set
$F_2(\La_{\partial G} P)$ is contained in the limit set
$\La_YQ=\{q\}$. So the map $f=F_2F_1^{-1}$ is well-defined on the
set of parabolic points of $\La_XG.$
 Furthermore   $f$ is $1{-}{\rm
to}{-}1$ at every conical point of $\La_X G$ \cite[Proposition
7.5.2]{Ge2}. Since all spaces are compacta the map $f$ is
continuous. It is also equivariant as $F_i$ are equivariant. So
the map $f$ satisfies the claim in this case.

The statement of d)  with the additional assumptions that  $G$ is
countable and $X$ and $Y$ are metrisable   was proved in
\cite{MOY}. Their proof uses the condition e) which was assumed to
be known in this case.

The statement e)  generalizes the last part of the main result of
\cite{Ya} to the case of non-finitely generated groups. It follows
from e) that for a $32$-action of $G$ on $X$ whose set of the
parabolic  points is $\Par_X$, there exists an equivariant
homeomorphism from $X$ to the Bowditch's boundary of the graph
$\G$ whose vertex set is $G\cup\Par_X.$}

\end{rems}

%
%
%
%
%
%
%
%
%

\end{document}